\documentclass[12pt]{article}
\usepackage{amsmath}
\usepackage{amssymb}
\allowdisplaybreaks

\def\q{\quad}
\def\qq{\qquad}
\def\qtq#1{\q\t{#1}\q}
\def\mod#1{\ (\text{\rm mod}\ #1)}
\def\t{\text}
\def\f{\frac}
\def\e{\equiv}
\def\b{\binom}
\def\sls#1#2{(\f{#1}{#2})}
 \def\ls#1#2{\big(\f{#1}{#2}\big)}
\def\Ls#1#2{\Big(\f{#1}{#2}\Big)}
\def\ap{\langle a\rangle_p}
\def\ep{\varepsilon_p}

\begin{document}
\leftline{}
\vskip0.2cm
 \centerline {\bf
Congruences involving the products of three binomial coefficients II}
\par\q\newline
\centerline{Zhi-Hong Sun}\newline \centerline{School of Mathematics
and Statistics} \centerline{Huaiyin Normal University}
\centerline{Huaian, Jiangsu 223300, P.R. China} \centerline{Email:
zhsun@hytc.edu.cn} \centerline{URL:
https://maths.hytc.edu.cn/szh.htm}
\vskip0.4cm
 \par{\bf Abstract} Let $p$ be an odd prime, $a,m\in\Bbb Z_{(p)}$ and $a(a+1)m\not\equiv 0\pmod p$. In this paper, we deduce the congruence for
 $(\sum_{k=0}^{p-1}\frac{\binom{2k}k\binom ak\binom{-1-a}k}{m^k})
 (\sum_{k=0}^{p-2}\frac{\binom{2k}k\binom ak\binom{-1-a}k}{m^k(k+1)})$ modulo $p^3$.
 As applications we obtain many congruences for
 \begin{align*}&\sum_{k=0}^{(p-1)/2}\frac{\binom{2k}k^3}
 {m^k(k+1)^r},\q
\sum_{k=0}^{p-2}\frac{\binom{2k}k^2\binom{4k}{2k}}{m^k(k+1)^r},
\q
\sum_{k=0}^{p-2}\frac{\binom{2k}k\binom{3k}k\binom{6k}{3k}}
{m^k(k+1)^r}\q(r=1,2,3)
\\&\sum_{k=0}^{(p-1)/2}k^j\frac{\binom{2k}k^3}
 {m^k},\q
\sum_{k=0}^{p-1}k^j\frac{\binom{2k}k^2\binom{4k}{2k}}{m^k},
\q
\sum_{k=0}^{p-1}k^j\frac{\binom{2k}k\binom{3k}k\binom{6k}{3k}}
{m^k}\q(j=2,3)
\end{align*} modulo $p^3$, which confirm the author's conjectures and improve previous results modulo $p^2$.

  \par\q

\section*{1. Introduction}
\par\q For a given real number $a$, the generalized binomial coefficient $\b ak$ is given by
$$\b a0=1\qtq{and}\b
ak=\f{a(a-1)\cdots(a-k+1)}{k!}\qtq{for}k=1,2,3,\ldots$$ One can easily check that $\b{-1/2}k=\b{2k}k\f 1{(-4)^k}$ and $p\mid \b {2k}k$ for
any odd prime $p$ satisfying $\f p2<k<p$. It is known (see for example [10-13] and [18, Lemma 21.3]) that
$$\b {-\f 13}k\b{-\f
23}k=\f{\b{2k}k\b{3k}k}{27^k}, \ \b {-\f 14}k\b{-\f
34}k=\f{\b{2k}k\b{4k}{2k}}{64^k}, \ \b {-\f 16}k\b{-\f
56}k=\f{\b{3k}{k}\b{6k}{3k}}{432^k}. \eqno{(1.1)}$$
\par
 For positive integers $a,b$ and $n$, the representation
 $n=ax^2+by^2$ means that
$n$ is represented by $ax^2+by^2(x,y\in\Bbb Z)$. Let $p$ be an odd prime. It is well known (see [1] and [3]) that
\begin{align*} & p=x^2+4y^2\q\t{for $p\e 1\mod 4$,}
\\& p=x^2+3y^2\q\t{for $p\e 1\mod 3$,}
\\&p=x^2+2y^2\q\t{for $p\e 1,3\mod 8$,}
\\&p=x^2+5y^2\q\t{for $p\e 1,9\mod {20}$,}
\\&p=x^2+7y^2\q\t{for $p\e 1,2,4\mod 7$.}
\end{align*}

\par Let $p>3$ be a prime and $m\in\Bbb Z_{(p)}$, where $Z_{(p)}$ is the set of those rational numbers whose denominators are not divisible by $p$. Based on Beukers' method [2], Sun and Ye [20] proved the congruences for
       $$\sum_{k=0}^{(p-1)/2}\b{2k}k^3\f 1{m^k},
       \q\sum_{k=0}^{p-1}\f{\b{2k}k^2\b{4k}{2k}}{m^k},
       \q\sum_{k=0}^{p-1}\f{\b{2k}k\b{3k}k\b{6k}{3k}}{m^k}$$     modulo $p^3$, where $p$ is an odd  prime represented by a certain binary quadratic form. For example, they showed that for any odd prime $p\e 1,2,4\mod 7$ and so $p=x^2+7y^2$,
       $$\sum_{k=0}^{(p-1)/2}\b{2k}k^3\e 4x^2-2p-\f {p^2}{4x^2}\mod {p^3},$$
       which was conjectured by the author in [14].
       \par Let $p>3$ be a prime and $m\in\Bbb Z_{(p)}$ with $m\not\e 0\mod p$. In [15,16,17,19], the author conejctured numerous congruences for
       \begin{align*}&\sum_{k=0}^{(p-1)/2}\frac{\binom{2k}k^3}
 {m^k(k+1)^r},\q \sum_{k=0}^{p-2}\frac{\binom{2k}k^2\binom{3k}{k}}{m^k(k+1)^r},
\q
\sum_{k=0}^{p-2}\frac{\binom{2k}k^2\binom{4k}{2k}}{m^k(k+1)^r},
\q
\sum_{k=0}^{p-2}\frac{\binom{2k}k\binom{3k}k\binom{6k}{3k}}
{m^k(k+1)^r}\q(r=1,2,3),
\\&\sum_{k=0}^{(p-1)/2}k^j\frac{\binom{2k}k^3}
 {m^k},\q
 \sum_{k=0}^{p-1}k^j\frac{\binom{2k}k^2\binom{3k}{k}}{m^k},
\q
\sum_{k=0}^{p-1}k^j\frac{\binom{2k}k^2\binom{4k}{2k}}{m^k},
\q
\sum_{k=0}^{p-1}k^j\frac{\binom{2k}k\binom{3k}k\binom{6k}{3k}}
{m^k}\q(j=2,3),
\\&\sum_{k=0}^{p-1}\frac{\binom{2k}k^3}
 {m^k(2k-1)}\q\t{and}\q \sum_{k=0}^{p-1}\frac{\binom{2k}k^3}
 {m^k(2k-1)^2}
\end{align*} modulo $p^3$. In [17], the author proved that
\begin{align*}
&\sum_{k=0}^{(p-1)/2}\f{\b{2k}k^3}{64^k(k+1)}\e 4x^2-2p
\mod {p^3}\q\t{for $p=x^2+4y^2\e 1\mod 4$,}
\\&\sum_{k=0}^{p-2}\f{\b{2k}k^2\b{3k}k}{108^k(k+1)}\e 4x^2-2p
\mod {p^3}\q\t{for $p=x^2+3y^2\e 1\mod 3$,}
\\&\sum_{k=0}^{p-2}\f{\b{2k}k^2\b{4k}{2k}}{256^k(k+1)}
\e 4x^2-2p\mod {p^3}\q\t{for $p=x^2+2y^2\e 1,3\mod 8$.}
\end{align*}
The corresponding congruences modulo $p^2$ were previously proved by the author's brother Z.W. Sun [22].
In [19], the author proved that
$$\Ls {-3}p\sum_{k=0}^{p-2}\f{\b{2k}k\b{3k}k\b{6k}{3k}}{12^{3k}(k+1)}\e
4x^2-2p\mod {p^3}\q\t{for $p=x^2+4y^2\e 1\mod 4$}$$
and deduced the congruences for
       \begin{align*}&\sum_{k=0}^{(p-1)/2}\frac{\binom{2k}k^3}
 {m^k(k+1)^r}\q (r=1,2,3),
 \ \sum_{k=0}^{(p-1)/2}k^2\frac{\binom{2k}k^3}
 {m^k},\\&  \sum_{k=0}^{(p-1)/2}k^3\frac{\binom{2k}k^3}
 {m^k},\ \sum_{k=0}^{p-1}\frac{\binom{2k}k^3}
 {m^k(2k-1)},\ \sum_{k=0}^{p-1}\frac{\binom{2k}k^3}
 {m^k(2k-1)^2}
 \end{align*}
    modulo $p^2$, where $\sls ap$ is the Legendre symbol.
    \par Let $p$ be an odd prime, $a,m\in\Bbb Z_{(p)}$ and $a(a+1)m\not\equiv 0\pmod p$. In this paper, we establish the congruence
    \begin{align*}&\Big(\sum_{k=0}^{p-1}\f{\b{2k}k\b ak\b{-1-a}k}{m^k}\Big)
 \Big(\sum_{k=0}^{p-2}\f{\b{2k}k\b ak\b{-1-a}k}{m^k(k+1)}\Big)
 \\&\e \Big(\sum_{k=0}^{p-1}\f{\b{2k}k\b ak\b{-1-a}k}{m^k}\Big)^2+\f{m-4}{4a(a+1)}
 \Big(\sum_{k=0}^{p-1}\f{k\b{2k}k\b ak\b{-1-a}k}{m^k}\Big)^2
 +p^2\f{a'(a'+1)}{a(a+1)}\mod {p^3},
 \end{align*}
 where $a'=(a-\ap)/p$ and $\ap$ is the least nonnegative residue of $a$ modulo $p$.
 As applications, we confirm many conjectures on
 \begin{align*}&\sum_{k=0}^{(p-1)/2}\frac{\binom{2k}k^3}
 {m^k(k+1)^r},\q
\sum_{k=0}^{p-2}\frac{\binom{2k}k^2\binom{4k}{2k}}{m^k(k+1)^r},
\q
\sum_{k=0}^{p-2}\frac{\binom{2k}k\binom{3k}k\binom{6k}{3k}}
{m^k(k+1)^r}\q(r=1,2,3),
\\&\sum_{k=0}^{(p-1)/2}k^j\frac{\binom{2k}k^3}
 {m^k},\q
\sum_{k=0}^{p-1}k^j\frac{\binom{2k}k^2\binom{4k}{2k}}{m^k},
\q
\sum_{k=0}^{p-1}k^j\frac{\binom{2k}k\binom{3k}k\binom{6k}{3k}}
{m^k}\q(j=2,3)
\end{align*} modulo $p^3$. In particular, we prove the following typical congruences:
\begin{align*}&\sum_{k=0}^{(p-1)/2}k^3\b{2k}k^3\e -\f{5408x^2}{27783}+\f{2992p}{9261}-\f{1774p^2}{27783x^2} \mod {p^3}\\&\qq\qq\qq\qq\qq\q\t{for $p\e 1,2,4\mod 7$ and so $p=x^2+7y^2$,}
\\&\sum_{k=0}^{(p-1)/2}\f{\b{2k}k^3}{k+1}
\e -44y^2+2p\mod {p^3}\q\t{for $p=x^2+7y^2\e 1,2,4\mod 7$,}
\\&(-1)^{\f{p-1}2}\sum_{k=0}^{p-1}\f{\b{2k}k^3}{(-64)^k(2k-1)}
\e p-3x^2\mod {p^3}\q\t{for $p=x^2+2y^2\e 1,3\mod 8$,}
\\&(-1)^{\f{p-1}2}\sum_{k=0}^{(p-1)/2}\f{\b{2k}k^3}{256^k(k+1)^2}\e -48y^2+8p\mod {p^3}\q\t{for $p=x^2+3y^2\e 1\mod 3$,}
\\&(-1)^{\f{p-1}4}
\sum_{k=0}^{(p-1)/2}\f{\b{2k}k^3}{(-512)^k(k+1)}
\e -32y^2+2p\pmod {p^3}\q\t{for $p=x^2+4y^2\e 1\mod 4$,}
\\&\sum_{k=0}^{p-2}\f{\b{2k}k^2\b{4k}{2k}}{(-1024)^k(k+1)}
\e\begin{cases} -32y^2+2p\pmod {p^3}&\t{if $p\e 1,9\pmod{20}$
and $p=x^2+5y^2$,}
\\ 16y^2-2p\pmod {p^3}&\t{if $p\e
3,7\pmod{20}$ and $2p=x^2+5y^2$,}\end{cases}
\\&\Ls{-2}p\sum_{k=0}^{p-2}\f{\b{2k}k\b{3k}k\b{6k}{3k}}
{(-32)^{3k}(k+1)}\e -26y^2+2p\mod {p^3}
\\&\qq\qq\qq\qq\qq\qq\t{for $p\e 1,3,4,5,9\mod {11}$ and so $4p=x^2+11y^2\ (x,y\in\Bbb Z)$.}
\end{align*}

\section*{2. A congruence for $\big(\sum_{k=0}^{p-1}\f{\b{2k}k\b ak\b{-1-a}k}{m^k}\big)
 \big(\sum_{k=0}^{p-2}\f{\b{2k}k\b ak\b{-1-a}k}{m^k(k+1)}\big)$
 modulo $p^3$}
\par\q\
Let
$(a)_k=a(a+1)\cdots(a+k-1)=(-1)^kk!\b{-a}k$ and
$${}_rF_s\biggl(\begin{matrix} a_1,\ldots,a_r\\
b_1,\ldots,b_s\end{matrix}
\biggm|z\biggr)=1+\sum_{k=1}^{\infty}\f{(a_1)_k\cdots(a_r)_k}
{(b_1)_k\cdots(b_s)_k}\cdot \f{z^k}{k!}.$$
Clausen's identity ([9, Eq.(2.5.7)]) states that
$${}_2F_1\biggl(\begin{matrix} a, \, b \\ a+b+\frac12 \end{matrix}
\biggm|z\biggr)^2 ={}_3F_2\biggl(\begin{matrix} 2a, \, 2b, \, a+b \\
a+b+\frac12, \, 2a+2b \end{matrix} \biggm|z\biggr).$$

\par\q\ {\bf Lemma 2.1} {\sl For $a\in\Bbb R$ let
$$f(z)=\sum_{k=0}^{\infty}\f{\b{2k}k\b ak\b{-1-a}k}{4^k}z^k
\q\t{and}\q g(z)=\sum_{k=0}^{\infty}\f{\b{2k}k\b ak\b{-1-a}k}{4^k(k+1)}z^k.$$
Then
$$a(a+1)(f(z)g(z)-f(z)^2)=\f{1-z}z(zf'(z))^2.$$}
\par{\it Proof.} Set
$$u(z)=\sum_{k=0}^{\infty}\b{\f a2}k\b{-\f{a+1}2}kz^k.$$
By Clausen's identity,
$$u(z)^2={}_2F_1\biggl(\begin{matrix} -\f a2, \, \f{a+1}2 \\ 1 \end{matrix}
\biggm|z\biggr)^2
={}_3F_2\biggl(\begin{matrix} -a, \, a+1, \, \f 12 \\
1, \, 1 \end{matrix} \biggm|z\biggr)=f(z).$$
Since $\b x{k+1}=\f{x-k}{k+1}\b xk$ we see that
\begin{align*}(k+1)^2\b {\f a2}{k+1}\b{-\f{a+1}2}{k+1}
&=\Big(\f a2-k\Big)\Big(-\f{a+1}2-k\Big)\b{\f a2}k\b{-\f{a+1}2}k
\\&=\Big(k(k-1)+\f 32k-\f{a(a+1)}4\Big)
\b{\f a2}k\b{-\f{a+1}2}k.\end{align*}
Let $u'(z)$ and $u''(z)$ be the derivative and second-order derivative of $u(z)$, respectively. Then
$$u'(z)=\sum_{k=1}^{\infty}k\b{\f a2}k\b{-\f{a+1}2}kz^{k-1}\q\t{and}
\q u''(z)=\sum_{k=2}^{\infty}k(k-1)\b{\f a2}k\b{-\f{a+1}2}kz^{k-2}.$$
Hence we have
\begin{align*}&z(1-z)u''(z)+\Big(1-\f 32z\Big)u'(z)+\f{a(a+1)}4u(z)
\\&=\sum_{k=1}^{\infty}k(k-1)\b{\f a2}k\b{-\f{a+1}2}k(z^{k-1}-z^k)+\sum_{k=1}^{\infty}k\b{\f a2}k\b{-\f{a+1}2}k(z^{k-1}-\f 32z^k)
\\&\q+\f{a(a+1)}4\sum_{k=0}^{\infty}\b{\f a2}k\b{-\f{a+1}2}kz^k
\\&=\sum_{k=0}^{\infty}\Big((k+1)k\b{\f a2}{k+1}\b{-\f{a+1}2}{k+1}
-k(k-1)\b{\f a2}k\b{-\f{a+1}2}k
+(k+1)\b{\f a2}{k+1}\b{-\f{a+1}2}{k+1}
\\&\q-\f 32k\b{\f a2}k\b{-\f{a+1}2}k+\f{a(a+1)}4\b{\f a2}k\b{-\f{a+1}2}k\Big)z^k
\\&=0.\end{align*}
Let $v(z)=a(a+1)u(z)^2+4z(1-z)u'(z)^2$. From the above we deduce that
\begin{align*}(zv(z))'&=a(a+1)u(z)^2+a(a+1)z\cdot 2u(z)u'(z)
\\&\q+(8z-12z^2)u'(z)^2+4z^2(1-z)\cdot 2u'(z)u''(z)
\\&=a(a+1)u(z)^2+8zu'(z)\Big(\f{a(a+1)}4u(z)+\Big(1-\f 32z\Big)u'(z)+z(1-z)u''(z)\Big)
\\&=a(a+1)u(z)^2=a(a+1)f(z).\end{align*}
Note that $(zg(z))'=f(z)$. We then obtain
$(zv(z)-a(a+1)zg(z))'=0$ and so $v(z)=a(a+1)g(z)$.
Hence,
\begin{align*}&a(a+1)f(z)g(z)-a(a+1)f(z)^2
\\&=a(a+1)f(z)g(z)-a(a+1)f(z)u(z)^2=4z(1-z)u(z)^2u'(z)^2
\\&=z(1-z)f'(z)^2=\f{1-z}z(zf'(z))^2.
\end{align*}
This proves the lemma.
\vskip0.2cm
\par{\bf Theorem 2.1} {\sl Let $p$ be an odd prime, $a,m\in\Bbb Z_{(p)}$ and $a(a+1)m\not\e 0\mod p$. Then
 \begin{align*}&\Big(\sum_{k=0}^{p-1}\f{\b{2k}k\b ak\b{-1-a}k}{m^k}\Big)
 \Big(\sum_{k=0}^{p-2}\f{\b{2k}k\b ak\b{-1-a}k}{m^k(k+1)}\Big)
 \\&\e \Big(\sum_{k=0}^{p-1}\f{\b{2k}k\b ak\b{-1-a}k}{m^k}\Big)^2+\f{m-4}{4a(a+1)}
 \Big(\sum_{k=0}^{p-1}\f{k\b{2k}k\b ak\b{-1-a}k}{m^k}\Big)^2
 +p^2\f{a'(a'+1)}{a(a+1)}\mod {p^3},
 \end{align*}
 where $a'=(a-\ap)/p$ and $\ap$ is the least nonnegative residue of $a$ modulo $p$.  }
 \vskip0.2cm
\par
{\it Proof}.
Set
$$F(z)=a(a+1)(f(z)g(z)-f(z)^2)-\f{1-z}z(zf'(z))^2\q \t{and}\q
c_k=\f{\b{2k}k\b ak\b{-1-a}k}{4^k}.$$
By Lemma 2.1,
$$F(z)=a(a+1)\Big(\Big(\sum_{k=0}^{\infty}c_kz^k\Big)
\Big(\sum_{k=0}^{\infty}\f{c_k}{k+1} z^k\Big)-\Big(\sum_{k=0}^{\infty}c_kz^k\Big)^2\Big)
-\f{1-z}z\Big(\sum_{k=0}^{\infty}kc_kz^k\Big)^2=0.$$
Let
$$f_p(z)=\sum_{k=0}^{p-1}c_kz^k,\q
g_p(z)=\sum_{k=0}^{p-2}\f{c_k}{k+1}z^k,
\q D_p(z)=\sum_{k=0}^{p-1}kc_kz^k$$
and
$$F_p(z)=a(a+1)(f_p(z)g_p(z)-f_p(z)^2)-\f{1-z}zD_p(z)^2.$$
Then $F_p(z)$ is a polynomial of $z$ with the degree less than $2p-1$. Let $[z^n]h(z)$ be the coefficient of $z^n$ in the power series expansion of $h(z)$.
For $n\in\{0,1,\ldots,p-2\}$ we see that
$$[z^n]F_p(z)=[z^n]\Big(a(a+1)(f(z)g(z)-f(z)^2)
-\f{1-z}z(zf'(z))^2\Big)
=[z^n]F(z)=0.$$
Also,
\begin{align*}[z^{p-1}]F_p(z)&=[z^{p-1}]\Big(F(z)-a(a+1)\f{c_{p-1}}pz^{p-1}\Big)
=[z^{p-1}]F(z)-a(a+1)\f{c_{p-1}}p
\\&=-a(a+1)\f{c_{p-1}}p.\end{align*}
For $0\le k\le p-1$ we have
$$\b ak\e \b{\ap}k\mod p\q\t{and}\q\b{-1-a}k\e \b{p-1-\ap}k\mod p.$$
 Also, $p\mid \b{2k}k$ for $\f p2<k<p$. Thus,
  $$c_k=\f{\b{2k}k\b ak\b{-1-a}k}{4^k}
  \e\begin{cases} 0\mod p&\t{if $\min\{\ap,p-1-\ap\}<k\le
  \f{p-1}2$,}
  \\ 0\mod {p^2}&\t{if $\f{p-1}2<k\le \max\{\ap,p-1-\ap\}$,}
  \\0\mod {p^3}&\t{if $\max\{\ap,p-1-\ap\}<k\le
   p-1$.}
  \end{cases}\eqno{(2.1)}$$
Suppose $i,j\in\{0,1,\ldots,p-1\}$ and $i+j\ge p$. If $i\le
\min\{\ap,p-1-\ap\}$, then $$j\ge p-i>p-1-\min\{\ap,p-1-\ap\}
=\max\{\ap,p-1-\ap\}$$ and so $c_j\e 0\mod {p^3}$ by (2.1). If $\min\{\ap,p-1-\ap\}<i\le \f{p-1}2$, then $c_i\e 0\mod p$, $j\ge p-i\ge \f{p+1}2$ and so
$c_j\e 0\mod {p^2}$ by (2.1). If $\f{p-1}2<i\le \max\{\ap,p-1-\ap\}$, then $c_i\e 0\mod {p^2}$, $j\ge p-i\ge p-\max\{\ap,p-1-\ap\}>\min\{\ap,p-1-\ap\}$ and so $c_j\e 0\mod p$ by (2.1). If $\max\{\ap,p-1-\ap\}<i\le p-1$, then $c_i\e 0\mod {p^3}$ by (2.1). Therefore, we always have
$$c_ic_j\e 0\mod {p^3}\q\t{for $i,j\in\{0,1,\ldots,p-1\}$ and $i+j\ge p$}.$$
It then follows that
$$[z^n]F_p(z)\e 0\mod {p^3}\q\t{for $n\in\{p,p+1,\ldots,2p-2\}$.}$$
Combining all the above yields
$$F_p(z)\e -a(a+1)\f{c_{p-1}}pz^{p-1}\mod {p^3}.$$
By [17, Lemma 3.2],
$$\f{c_{p-1}}p\e -\f{a'(a'+1)}{a(a+1)}p^2\mod {p^3}.$$
Hence, for $z\in\Bbb Z_{(p)}$ with $z\not\e 0\mod p$,
$$F_p(z)\e -a(a+1)\f{c_{p-1}}pz^{p-1}\e a'(a'+1)p^2\mod {p^3}.$$
Taking $z=\f 4m$ in the congruence yields the result.
\vskip0.2cm
\par\q{\bf Corollary 2.1} {\sl Let $p$ be an odd prime, $m\in\Bbb Z_{(p)}$ and $a(a+1)m\not\e 0\mod p$. Then
 \begin{align*}&\Big(\sum_{k=0}^{(p-1)/2}
 \f{\b{2k}k^3}{(16m)^k}\Big)
 \Big(\sum_{k=0}^{(p-1)/2}\f{\b{2k}k^3}{(16m)^k(k+1)}\Big)
 \\&\e \Big(\sum_{k=0}^{(p-1)/2}\f{\b{2k}k^3}{(16m)^k}\Big)^2
 +(4-m) \Big(\sum_{k=0}^{(p-1)/2}\f{k\b{2k}k^3}{(16m)^k}\Big)^2
 +p^2\mod {p^3}.
 \end{align*}}
\par{\it Proof.} Since $p\mid \b{2k}k$ for $\f p2<k<p$, taking $a=-\f 12$ in Theorem 2.1 and noting that $\b{-\f 12}k=\b{2k}k\f 1{(-4)^k}$ gives the result.
\vskip0.2cm
\par{\bf Corollary 2.2} {\sl Let $p$ be a prime greater than $3$, $m\in\Bbb Z_{(p)}$ and $m\not\e 0\mod p$. Then
 \begin{align*}&\Big(\sum_{k=0}^{p-1}\f{\b{2k}k^2\b{3k}k}
 {(27m)^k}\Big)
 \Big(\sum_{k=0}^{p-2}\f{\b{2k}k^2\b{3k}k}{(27m)^k(k+1)}\Big)
 \\&\e \Big(\sum_{k=0}^{p-1}\f{\b{2k}k^2\b{3k}k}{(27m)^k}\Big)^2
 +\f{9(4-m)}{8}
 \Big(\sum_{k=0}^{p-1}\f{k\b{2k}k^2\b{3k}k}{(27m)^k}\Big)^2
 +p^2\mod {p^3}.
 \end{align*}}
 \par{\it Proof.} Taking $a=-\f 13$ in Theorem 2.1 and then applying (1.1) gives the result.
 \vskip 0.2cm
 \vskip0.2cm
\par{\bf Corollary 2.3} {\sl Let $p>3$ be a prime, $m\in\Bbb Z_{(p)}$ and $m\not\e 0\mod p$. Then
 \begin{align*}&\Big(\sum_{k=0}^{p-1}\f{\b{2k}k^2\b{4k}{2k} }{(64m)^k}\Big)
 \Big(\sum_{k=0}^{p-2}\f{\b{2k}k^2\b{4k}{2k}}{(64m)^k(k+1)}\Big)
 \\&\e \Big(\sum_{k=0}^{p-1}\f{\b{2k}k^2\b{4k}{2k}}{(64m)^k}\Big)^2
 +\f{4(4-m)}{3}
 \Big(\sum_{k=0}^{p-1}\f{k\b{2k}k^2\b{4k}{2k}}{(64m)^k}\Big)^2
 +p^2\mod {p^3}.
 \end{align*}}
 \par{\it Proof.} Taking $a=-\f 14$ in Theorem 2.1 and then applying (1.1) yields the result.
 \vskip0.2cm
\par{\bf Corollary 2.4} {\sl Let $p>5$ be a prime, $m\in\Bbb Z_{(p)}$ and $m\not\e 0\mod p$. Then
 \begin{align*}&\Big(\sum_{k=0}^{p-1}\f{\b{2k}k\b{3k}k\b{6k}{3k}}
 {(432m)^k}\Big)
 \Big(\sum_{k=0}^{p-2}\f{\b{2k}k\b {3k}k\b{6k}{3k}}{(432m)^k(k+1)}\Big)
 \\&\e \Big(\sum_{k=0}^{p-1}\f{\b{2k}k\b{3k}k\b{6k}{3k}
 }{(432m)^k}\Big)^2+\f{9(4-m)}{5}
 \Big(\sum_{k=0}^{p-1}\f{k\b{2k}k\b{3k}k\b{6k}{3k}}{(432m)^k}\Big)^2
 +p^2\mod {p^3}.
 \end{align*}}
 \par{\it Proof.} Taking $a=-\f 16$ in Theorem 2.1 and then applying (1.1) yields the result.
 \vskip0.2cm

 \par{\bf Lemma 2.2 ([19, Theorem 4.1])} {\sl Let $p$ be an odd prime, $a,m\in\Bbb Z_{(p)}$,
$a\not\e 0,-1\mod p$ and $m\not\e 0\mod p$. Then
\begin{align*} (m-4)&\sum_{k=0}^{p-1}\f{k^2\b ak\b{-1-a}k\b{2k}k}{m^k}
\e 2\sum_{k=0}^{p-1}\f{k\b ak\b{-1-a}k\b{2k}k}{m^k}
-4a(a+1)\sum_{k=0}^{p-1}\f{\b ak\b{-1-a}k\b{2k}k}{m^k}
\\&\qq\qq\qq\qq\qq+2a(a+1)\sum_{k=0}^{p-2} \f{\b ak\b{-1-a}k\b{2k}k}
{m^k(k+1)}\mod {p^3}\end{align*} and
\begin{align*}&(m-4)^2\sum_{k=0}^{p-1}\f{k^3\b ak\b{-1-a}k\b{2k}k}{m^k}
\e ((2-4a(a+1))(m-4)+12)\sum_{k=0}^{p-1}\f{k\b
ak\b{-1-a}k\b{2k}k}{m^k}\\&\q -2a(a+1)(m+8)\sum_{k=0}^{p-1}\f{\b
ak\b{-1-a}k\b{2k}k}{m^k}+12a(a+1)\sum_{k=0}^{p-2}\f{\b
ak\b{-1-a}k\b{2k}k}{m^k(k+1)}\mod{p^3}.
\end{align*}}
 \par{\bf Lemma 2.3 ([19, Theorem 3.1])} {\sl Let $p$ be an odd prime, $a,m\in\Bbb Z_{(p)}$, $m\not\e 0\mod p$
  and
 $a\not \e 0,-1\mod p$. Then
\begin{align*} 2a(a+1)\sum_{k=0}^{p-2}\f {\b ak\b{-1-a}k\b{2k}k}{m^k(k+1)^2}
&\e (m-4)\sum_{k=0}^{p-1}\f {k\b ak\b{-1-a}k\b{2k}k}{m^k}
+2\sum_{k=0}^{p-1}\f {\b ak\b{-1-a}k\b{2k}k}{m^k} \\&\q+(4a(a+1)-2)
\sum_{k=0}^{p-2}\f {\b
ak\b{-1-a}k\b{2k}k}{m^k(k+1)}\mod{p^3}\end{align*} and
\begin{align*} &2a(a+1)\sum_{k=0}^{p-2}\f {\b ak\b{-1-a}k\b{2k}k}{m^k(k+1)^3}
\\&\e-m+\Big(2-\f 1{a(a+1)}\Big)(m-4) \sum_{k=0}^{p-1}\f {k\b
ak\b{-1-a}k\b{2k}k}{m^k} +\Big(m-\f
2{a(a+1)}\Big)\sum_{k=0}^{p-1}\f {\b ak\b{-1-a}k\b{2k}k}{m^k}
\\&\q+\Big(8a(a+1)-2+\f 2{a(a+1)}\Big) \sum_{k=0}^{p-2}\f {\b
ak\b{-1-a}k\b{2k}k}{m^k(k+1)} \mod {p^3}.
\end{align*}}

\par{\bf Lemma 2.4 ([19,Theorem 5.1])} {\sl Suppose that $p$ is an odd prime, $a,m\in\Bbb
Z_{(p)}$, $m\not\e 0\mod p$ and $a\not\e 0,-1\mod p$. Then}
\begin{align*}\sum_{k=0}^{p-1}\f{\b ak\b{-1-a}k\b{2k}k}{m^k(2k-1)}
&\e \Big(\f 8m-2\Big)\sum_{k=0}^{p-1}\f{k\b ak\b{-1-a}k\b{2k}k}{m^k}
-\sum_{k=0}^{p-1}\f{\b ak\b{-1-a}k\b{2k}k}{m^k}
\\&\q\ -\f{8a(a+1)}m\sum_{k=0}^{p-2}\f{\b ak\b{-1-a}k\b{2k}k}{m^k(k+1)}\mod {p^3}.
\end{align*}
\par{\bf Theorem 2.2} {\sl Let $p$ be an odd prime, $a,m\in\Bbb Z_{(p)}$ and $a(a+1)m(m-4)\not\e 0\mod p$. Let $a'=(a-\ap)/p$,
$$S_j=\sum_{k=0}^{p-1}k^j\f{\b{2k}k\b ak\b{-1-a}k}{m^k},\q T_j=\sum_{k=0}^{p-2}\f{\b{2k}k\b ak\b{-1-a}k}{m^k(k+1)^j}\q (j=0,1,2,3).$$
Then
\begin{align*}&(m-4)S_3\e -2a(a+1)S_0+(2-4a(a+1))S_1+6S_2\mod {p^3},\tag {2.2}
\\&(m-4)S_0S_2\e -2a(a+1)S_0^2+2S_0S_1+\f{m-4}2S_1^2+2a'(a'+1)p^2\mod {p^3},\tag{2.3}
\\&(m-4)^2S_0S_3\e -2a(a+1)(m+2)S_0^2+((2-4a(a+1))(m-4)+12)S_0S_1
\\&\q\qq\qq\qq\qq+3(m-4)S_1^2+12a'(a'+1)p^2\mod {p^3},\tag {2.4}
\\&2a(a+1)S_0T_2
\e (m-4)S_0S_1+4a(a+1)S_0^2+\Big(1-\f 1{2a(a+1)}\Big)(m-4)S_1^2
\\&\q\qq\qq\qq\q+\Big(4-\f 2{a(a+1)}\Big)a'(a'+1)p^2\mod {p^3},\tag{2.5}
\\&2a(a+1)S_0T_3\e  -mS_0+\Big(2-\f 1{a(a+1)}\Big)(m-4)S_0S_1 +\Big(8a(a+1)+m-2\Big)S_0^2
 \\&\q\qq\qq+\Big(8a(a+1)-2+\f 2{a(a+1)}\Big)\Big(\f{m-4}{4a(a+1)}S_1^2+\f{a'(a'+1)}{a(a+1)}p^2\Big)
 \mod {p^3}\tag {2.6}\end{align*}
 and
 \begin{align*}
 S_0\sum_{k=0}^{p-1}\f{\b{2k}k\b ak\b{-1-a}k}{m^k(2k-1)}&
 \e\Big(\f 8m-2\Big)S_0S_1-\Big(1+\f{8a(a+1)}m\Big)S_0^2
 \\&\q-\f{2(m-4)}mS_1^2-\f {8a'(a'+1)}mp^2\mod{p^3}.\tag{2.7}
\end{align*}
Hence, for $S_0\e 0\mod {p}$ we have
\begin{align*}&(m-4)^2S_0S_3\e 6(m-4)S_0S_2
\e 3(m-4)S_1^2+12a'(a'+1)p^2\mod {p^3},
\\&2a(a+1)S_0T_2
\e \Big(1-\f 1{2a(a+1)}\Big)(m-4)S_1^2
+\Big(4-\f 2{a(a+1)}\Big)a'(a'+1)p^2\mod {p^3},
\\&2a(a+1)S_0T_3\e  -mS_0+\Big(8a(a+1)-2+\f 2{a(a+1)}\Big)\Big(\f{m-4}{4a(a+1)}S_1^2+\f{a'(a'+1)}{a(a+1)}p^2\Big)
 \mod {p^3},
 \\&S_0\sum_{k=0}^{p-1}\f{\b{2k}k\b ak\b{-1-a}k}{m^k(2k-1)}
 \e -\f{2(m-4)}mS_1^2-\f {8a'(a'+1)}mp^2\mod{p^3}.
\end{align*}}
\par{\it Proof.}
By Lemma 2.2,
\begin{align*}
&(m-4)S_2\e 2S_1-4a(a+1)S_0+2a(a+1)T_1\mod {p^3},\tag{2.8}
\\&(m-4)^2S_3\e ((2-4a(a+1))(m-4)+12)S_1
\\&\qq\qq\qq-2a(a+1)(m+8)S_0+12a(a+1)T_1\mod {p^3}.\tag {2.9}
\end{align*}
Thus,
$$(m-4)^2S_3-6(m-4)S_2\e (2-4a(a+1))(m-4)S_1-2a(a+1)(m-4)S_0\mod {p^3}.$$
This yields (2.2).
\par From Theorem 2.1,
$$S_0T_1\e S_0^2+\f{m-4}{4a(a+1)}S_1^2+\f{a'(a'+1)}{a(a+1)}p^2\mod {p^3}.\eqno{(2.10)}$$
This together with (2.8) yields
\begin{align*}
(m-4)S_0S_2&\e 2S_0S_1-4a(a+1)S_0^2+2a(a+1)S_0T_1
\\&\e 2S_0S_1-4a(a+1)S_0^2+2a(a+1) S_0^2+\f{m-4}2S_1^2+2a'(a'+1)p^2\mod {p^3},
\end{align*}
which proves (2.3). Combining (2.9) with (2.10) yields
\begin{align*}
(m-4)^2S_0S_3&\e ((2-4a(a+1))(m-4)+12)S_0S_1
-2a(a+1)(m+8)S_0^2+12a(a+1)S_0T_1
\\&\e ((2-4a(a+1))(m-4)+12)S_0S_1
-2a(a+1)(m+8)S_0^2+12a(a+1)S_0^2
\\&\q+3(m-4)S_1^2+12a'(a'+1)p^2 \mod {p^3},
\end{align*}
which proves (2.4).
\par From Lemma 2.3,
\begin{align*} &2a(a+1)T_2
\e (m-4)S_1+2S_0+(4a(a+1)-2)T_1
\mod{p^3},
\\&2a(a+1)T_3\e -m+\Big(2-\f 1{a(a+1)}\Big)(m-4)S_1 +\Big(m-\f
2{a(a+1)}\Big)S_0
\\&\q\qq\qq\qq+\Big(8a(a+1)-2+\f 2{a(a+1)}\Big)T_1 \mod {p^3}.
\end{align*}
This together with (2.8) yields
\begin{align*}
2a(a+1)S_0T_2&\e
(m-4)S_0S_1+2S_0^2+(4a(a+1)-2)
\Big(S_0^2+\f{m-4}{4a(a+1)}S_1^2+\f{a'(a'+1)}{a(a+1)}p^2\Big)
\\&=(m-4)S_0S_1+4a(a+1)S_0^2+\Big(1-\f 1{2a(a+1)}\Big)(m-4)S_1^2
\\&\q+\Big(4-\f 2{a(a+1)}\Big)a'(a'+1)p^2\mod {p^3}.
\end{align*}
Similarly,
\begin{align*}
&2a(a+1)S_0T_3\\&\e -mS_0+\Big(2-\f 1{a(a+1)}\Big)(m-4)S_0S_1 +\Big(m-\f
2{a(a+1)}\Big)S_0^2
\\&\q+\Big(8a(a+1)-2+\f 2{a(a+1)}\Big)\Big(S_0^2+\f{m-4}{4a(a+1)}S_1^2+\f{a'(a'+1)}{a(a+1)}p^2\Big)
 \\&\e -mS_0+\Big(2-\f 1{a(a+1)}\Big)(m-4)S_0S_1 +\Big(8a(a+1)+m-2\Big)S_0^2
 \\&\q+\Big(8a(a+1)-2+\f 2{a(a+1)}\Big)\Big(\f{m-4}{4a(a+1)}S_1^2+\f{a'(a'+1)}{a(a+1)}p^2\Big)
 \mod {p^3}.
\end{align*}
By Lemma 2.4 and (2.10),
\begin{align*}
 &S_0\sum_{k=0}^{p-1}\f{\b{2k}k\b ak\b{-1-a}k}{m^k(2k-1)}
 \\&
 \e \Big(\f 8m-2\Big)S_0S_1-S_0^2-\f{8a(a+1)}m\Big(S_0^2+\f{m-4}{4a(a+1)}S_1^2
 +\f{a'(a'+1)}{a(a+1)}p^2\Big)
 \\&=\Big(\f 8m-2\Big)S_0S_1-\Big(1+\f{8a(a+1)}m\Big)S_0^2
 -\f{2(m-4)}mS_1^2-\f {8a'(a'+1)}mp^2\mod{p^3}.
\end{align*}
\par
For $S_0\e 0\mod p$ we have
 \begin{align*}&\sum_{k=0}^{p-1}\b ak\b{-1-a}k\Ls{-1\pm\sqrt{1-\f 4m}}2^k\e 0\mod p,
 \\& S_0\e 0\mod {p^2}\q{and}\q S_1\e 0\mod p\tag {2.11}
 \end{align*} by [19, Lemma 2.2 and Theorem 2.1] ($S_1\e 0\mod p$ also follows from (2.10)). Hence $(m-4)S_3\e 6S_2\mod p$ by (2.2). Since $S_0^2\e S_0S_1\e 0\mod {p^3}$, from (2.3)-(2.7) we deduce the remaining congruences.
 \vskip0.2cm
 \par{\bf Lemma 2.5 ([19, Theorem 5.2])} {\sl Let $p$ be an odd prime and $m\in\Bbb Z_{(p)}$ with $m\not\e 0\mod p$. Then
 \begin{align*}\sum_{k=0}^{p-1}\f{\b{2k}k^3}{(16m)^k(2k-1)^2}
&\e \Big(4-\f {16}m\Big)\sum_{k=0}^{(p-1)/2}
\f{k\b{2k}k^3}{(16m)^k}+\Big(1+\f 4m\Big)\sum_{k=0}^{(p-1)/2}
\f{\b{2k}k^3}{(16m)^k}\\&\q-\f 6m\sum_{k=0}^{(p-1)/2}
\f{\b{2k}k^3}{(16m)^k(k+1)}\mod {p^3}.\end{align*}}
 \par{\bf Theorem 2.3} {\sl Let $p$ be an odd prime and $m\in\Bbb Z_{(p)}$ with $m\not\e 0\mod p$. Then
 \begin{align*}
 S_0 \Big(\sum_{k=0}^{p-1}\f{\b{2k}k^3}{(16m)^k(2k-1)^2}\Big)
 \e
 \Big(1-\f 2m\Big)S_0^2+\Big(4-\f{16}m\Big)
 S_0S_1 +\f{6(m-4)}mS_1^2 -\f 6mp^2\mod {p^3},
  \end{align*}
  where
  $$S_0=\sum_{k=0}^{(p-1)/2}\f{\b{2k}k^3}{(16m)^k}
  \q{and}\q S_1=\sum_{k=0}^{(p-1)/2}\f{k\b{2k}k^3}{(16m)^k}.$$
 Hence, for $S_0\e 0\mod p$ we have}
 $$S_0\Big(\sum_{k=0}^{p-1}\f{\b{2k}k^3}{(16m)^k(2k-1)^2}\Big)
 \e \f{6(m-4)}mS_1^2
 -\f 6mp^2\mod {p^3}.$$
\par{\it Proof.} By Lemma 2.5 and Theorem 2.1 (with $a=-\f 12$),
\begin{align*}&S_0\Big(\sum_{k=0}^{p-1}\f{\b{2k}k^3}{(16m)^k(2k-1)^2}\Big)
\\&\e \Big(4-\f{16}m\Big)S_0S_1+\Big(1+\f 4m\Big)S_0^2
-\f 6m\big(S_0^2+(4-m)S_1^2+p^2\big)
\\&=\Big(1-\f 2m\Big)S_0^2+\Big(4-\f{16}m\Big)
 S_0S_1 +\f{6(m-4)}mS_1^2 -\f 6mp^2\mod {p^3}.
  \end{align*}
For $S_0\e 0\mod p$ we have $S_0\e 0\mod {p^2}$ and $S_1\e 0\mod p$ by (2.11). Thus, the remaining part follows.

 \section*{3. Congruences for $\sum_{k=0}^{(p-1)/2}\f{\b{2k}k^3}{m^k(k+1)}$ modulo $p^3$}
  \par\q\ {\bf Lemma 3.1 ([19, Corollary 3.1])} {\sl Let $p$ be an odd prime,
$m\in\Bbb Z_{(p)}$ and
 $m\not \e 0\mod p$. Then}
\begin{align*}&\sum_{k=0}^{\f{p-1}2}\f{\b{2k}k^3}{(16m)^k(k+1)^2}
\e (8-2m)\sum_{k=0}^{\f{p-1}2}\f{k\b{2k}k^3}{(16m)^k}
-4\sum_{k=0}^{\f{p-1}2}\f{\b{2k}k^3}{(16m)^k}
+6\sum_{k=0}^{\f{p-1}2}\f{\b{2k}k^3}{(16m)^k(k+1)}\mod {p^3},
\\&\sum_{k=0}^{(p-1)/2}\f{\b{2k}k^3}{(16m)^k(k+1)^3}
\e 2m-12(m-4)\sum_{k=0}^{(p-1)/2}\f{k\b{2k}k^3}{(16m)^k}
-2(m+8)\sum_{k=0}^{(p-1)/2}\f{\b{2k}k^3}{(16m)^k}
\\&\q\qq\qq\qq\qq\qq+24\sum_{k=0}^{(p-1)/2}\f{\b{2k}k^3}{(16m)^k(k+1)}\mod {p^3}.
\end{align*}

\par{\bf Lemma 3.2 ([19, Corollary 4.1])} {\sl Let $p$ be an odd
prime, $m\in\Bbb Z_{(p)}$ and $m\not\e 0\mod p$. Then
$$(m-4)\sum_{k=0}^{\f{p-1}2}\f{k^2\b{2k}k^3}{(16m)^k}
\e 2\sum_{k=0}^{\f{p-1}2}\f{k\b{2k}k^3}{(16m)^k}+
\sum_{k=0}^{\f{p-1}2}\f{\b{2k}k^3}{(16m)^k}-\f
12\sum_{k=0}^{\f{p-1}2}\f{\b{2k}k^3}{(16m)^k(k+1)}\mod {p^3}$$ and
\begin{align*} &(m-4)^2\sum_{k=0}^{(p-1)/2}\f{k^3\b{2k}k^3}{(16m)^k} \\&\e
3m\sum_{k=0}^{(p-1)/2}\f{k\b{2k}k^3}{(16m)^k}+\f {m+8}{2}
\sum_{k=0}^{(p-1)/2}\f{\b{2k}k^3}{(16m)^k}-3
\sum_{k=0}^{(p-1)/2}\f{\b{2k}k^3}{(16m)^k(k+1)}\mod
{p^3}.\end{align*}}

 \par{\bf Lemma 3.3 ([19, Corollary 5.1])}  {\sl Let $p$ be an odd
prime, $m\in\Bbb Z_{(p)}$ and $m\not\e 0\mod p$. Then}
\begin{align*}&\sum_{k=0}^{p-1}\f{\b{2k}k^3}{(16m)^k(2k-1)}
\\&\e \Big(\f 8m-2\Big)\sum_{k=0}^{(p-1)/2}
\f{k\b{2k}k^3}{(16m)^k}-\sum_{k=0}^{(p-1)/2} \f{\b{2k}k^3}{(16m)^k}
+\f{2}m\sum_{k=0}^{(p-1)/2} \f{\b{2k}k^3}{(16m)^k(k+1)} \mod
{p^3}.\end{align*}

 \par{\bf Theorem 3.1}  {\sl Let $p$ be an odd prime such that $p\e 1,2,4\mod 7$ and so $p=x^2+7y^2$ for some $x,y\in\Bbb Z$. Then
 \begin{align*}
 &\sum_{k=0}^{(p-1)/2}\f{\b{2k}k^3}{k+1}\e -44y^2+2p\mod {p^3},
 \\&\sum_{k=0}^{(p-1)/2}\f{\b{2k}k^3}{(k+1)^2}\e -68y^2+p-\f{p^2}{4y^2}\mod {p^3},
 \\&\sum_{k=0}^{(p-1)/2}\f{\b{2k}k^3}{(k+1)^3}\e \f 18-\f{201}2y^2-\f 94p-\f{39p^2}{32y^2}\mod {p^3},
 \\&\sum_{k=0}^{(p-1)/2}k^2\b{2k}k^3\e \f{736x^2}{1323}-\f{272p}{441}+\f{20p^2}{1323x^2}\mod {p^3},
 \\&\sum_{k=0}^{(p-1)/2}k^3\b{2k}k^3\e -\f{5408x^2}{27783}+\f{2992p}{9261}-\f{1774p^2}{27783x^2} \mod {p^3},
 \\&\sum_{k=0}^{p-1}\f{\b{2k}k^3}{2k-1} \e  -36y^2+14p-\f 7{4y^2}p^2\mod {p^3},
 \\&\sum_{k=0}^{p-1}\f{\b{2k}k^3}{(2k-1)^2}\e -284y^2+34p+\f {23}{4y^2}p^2\mod {p^3}
 \end{align*}
 and
 \begin{align*}
 &(-1)^{\f{p-1}2}\sum_{k=0}^{(p-1)/2}\f{\b{2k}k^3}{4096^k(k+1)}\e 72y^2+2p\mod {p^3},
 \\&(-1)^{\f{p-1}2}\sum_{k=0}^{(p-1)/2}\f{\b{2k}k^3}
 {4096^k(k+1)^2}\e -1136y^2+64p+\f{2p^2}{y^2}\mod {p^3},
 \\&(-1)^{\f{p-1}2}\sum_{k=0}^{(p-1)/2}\f{\b{2k}k^3}
 {4096^k(k+1)^3}\e 512(-1)^{\f{p-1}2}+6432y^2-648p-\f{6p^2}{y^2} \mod {p^3},
 \\&(-1)^{\f{p-1}2}\sum_{k=0}^{(p-1)/2}\f{k^2\b{2k}k^3}{4096^k}\e \f{43x^2}{1323}-\f{13}{441}p-\f{p^2}{1323x^2}\mod {p^3},
 \\&(-1)^{\f{p-1}2}\sum_{k=0}^{(p-1)/2}\f{k^3\b{2k}k^3}{4096^k}
 \e \f{169x^2}{55566}-\f{31}{74088}p-\f{71p^2}{444528x^2} \mod {p^3},
 \\&(-1)^{\f{p-1}2}\sum_{k=0}^{p-1}\f{\b{2k}k^3}{4096^k(2k-1)}
\e 22y^2-\f 74p-\f{7p^2}{256y^2}\mod {p^3},
\\&(-1)^{\f{p-1}2}\sum_{k=0}^{p-1}\f{\b{2k}k^3}{4096^k(2k-1)^2}
\e -17y^2+\f{97}{64}p+\f{5p^2}{256y^2}\mod{p^3}.
 \end{align*}}
 \par{\it Proof.} From [20, Theorem 1.1],
 $$\sum_{k=0}^{(p-1)/2}\b{2k}k^3\e 4x^2-2p-\f{p^2}{4x^2}\mod {p^3}.\eqno{(3.1)}$$
 In [21], Z.W. Sun showed that
 $$\sum_{k=0}^{(p-1)/2}(21k+8)\b{2k}k^3\e 8p\mod {p^3}.$$
 This together with (3.1) gives
 $$\sum_{k=0}^{(p-1)/2}k\b{2k}k^3\e
 \f 8{21}\Big(3p-4x^2+\f{p^2}{4x^2}\Big)
 \mod {p^3}.\eqno{(3.2)}$$
  Taking $m=\f 1{16}$ in Corollary 2.1 and then applying (3.1) and (3.2) gives
 \begin{align*}
 \sum_{k=0}^{(p-1)/2}\f{\b{2k}k^3}{k+1}
 &\e \sum_{k=0}^{(p-1)/2}\b{2k}k^3+
 \f{(4-\f 1{16})(\sum_{k=0}^{(p-1)/2}k\b{2k}k^3)^2+p^2}
 {\sum_{k=0}^{(p-1)/2}\b{2k}k^3}
 \\&\e 4x^2-2p-\f{p^2}{4x^2}+
 \f{\f{63}{16}\cdot\f{8^2}{21^2}(3p-4x^2+\f{p^2}{4x^2})^2+p^2}
 {4x^2-2p-\f{p^2}{4x^2}}
 \\&\e 4x^2-2p-\f{p^2}{4x^2}+\f{\f 47(16x^4-24x^2p+7p^2)+p^2}
 {4x^2-2p-\f{p^2}{4x^2}}
 \mod {p^3}.
 \end{align*}
It is easily seen that
\begin{align*}
&\Big(3p-4x^2+\f{p^2}{4x^2}\Big)^2\e 16x^4-24x^2p+7p^2\mod {p^3},
\\&\Big(4x^2-2p-\f{p^2}{4x^2}\Big)\Big(\f 1{4x^2}+\f p{8x^4}
+\f {5p^2}{64x^6}\Big)=1-\f{3p^3}{16x^6}
-\f{5p^4}{256x^8}\e 1\mod {p^3}.\tag{3.3}\end{align*}
Thus,
\begin{align*}
 \sum_{k=0}^{(p-1)/2}\f{\b{2k}k^3}{k+1}
 &\e 4x^2-2p-\f{p^2}{4x^2}+
 \Big(\f 1{4x^2}+\f p{8x^4}
+\f {5p^2}{64x^6}\Big)\Big(\f{4}{7}(16x^4-24x^2p+7p^2)+p^2\Big)
\\&\e 4x^2-2p-\f{p^2}{4x^2}+\f{16}7x^2-\f{16}7p+\f{p^2}{4x^2}
=\f{44x^2-30p}7
\\&=\f{44(p-7y^2)-30p}7=-44y^2+2p
\mod {p^3}.\tag{3.4}
\end{align*}
This together with (3.1), (3.2) and Lemma 3.1 (with $m=\f 1{16}$) yields
\begin{align*}
\sum_{k=0}^{\f{p-1}2}\f{\b{2k}k^3}{(k+1)^2}
&\e\Big(8-\f 2{16}\Big)\sum_{k=0}^{\f{p-1}2}k\b{2k}k^3
-4\sum_{k=0}^{\f{p-1}2}\b{2k}k^3+6\sum_{k=0}^{\f{p-1}2}\f{\b{2k}k^3}{k+1}
\\&\e \f{63}8\cdot \f 8{21}\Big(3p-4x^2+\f{p^2}{4x^2}\Big)
-4\Big(4x^2-2p-\f{p^2}{4x^2}\Big)+6(-44y^2+2p)
\\&\e -68y^2+p-\f{p^2}{4y^2}\mod {p^3}
\end{align*}
and
\begin{align*}
\sum_{k=0}^{\f{p-1}2}\f{\b{2k}k^3}{(k+1)^3}
&\e \f 2{16}-12
\Big(\f 1{16}-4\Big)\sum_{k=0}^{\f{p-1}2}k\b{2k}k^3
-2\Big(\f 1{16}+8\Big)\sum_{k=0}^{\f{p-1}2}\b{2k}k^3+24\sum_{k=0}^{\f{p-1}2}\f{\b{2k}k^3}{k+1}
\\&\e \f 18+\f{189}4\Big(3p-4x^2+\f{p^2}{4x^2}\Big)
-\f{129}8\Big(4x^2-2p-\f{p^2}{4x^2}\Big)+24(-44y^2+2p)
\\&\e \f 18-\f{201}2y^2-\f 94p-\f{39p^2}{32y^2}\mod {p^3}.
\end{align*}
Similarly, from (3.1), (3.2), (3.4) and Lemma 3.2 (with $m=\f 1{16}$) we derive
\begin{align*}
\sum_{k=0}^{\f{p-1}2}k^2\b{2k}k^3
&\e -\f{32}{63}\sum_{k=0}^{\f{p-1}2}k\b{2k}k^3-\f{16}{63}\sum_{k=0}^{\f{p-1}2}\b{2k}k^3
+\f 8{63}\sum_{k=0}^{\f{p-1}2}\f{\b{2k}k^3}{k+1}
\\&\e -\f{32}{63}\Big(3p-4x^2+\f{p^2}{4x^2}\Big)
-\f{16}{63}\Big(4x^2-2p-\f{p^2}{4x^2}\Big)+\f 8{63}(-44y^2+2p)
\\&\e  \f{736x^2}{1323}-\f{272p}{441}+\f{20p^2}{1323x^2}\mod{p^3}
\end{align*}
and
\begin{align*}
\sum_{k=0}^{\f{p-1}2}k^3\b{2k}k^3
&\e \f{16}{1323}\sum_{k=0}^{\f{p-1}2}k\b{2k}k^3
+\f{344}{1323}\sum_{k=0}^{\f{p-1}2}\b{2k}k^3
-\f {256}{1323}\sum_{k=0}^{\f{p-1}2}\f{\b{2k}k^3}{k+1}
\\&\e \f{16}{1323}\Big(3p-4x^2+\f{p^2}{4x^2}\Big)
+\f{344}{1323}\Big(4x^2-2p-\f{p^2}{4x^2}\Big)-\f {256}{1323}(-44y^2+2p)
\\&\e  -\f{5408x^2}{27783}+\f{2992p}{9261}-\f{1774p^2}{27783x^2} \mod {p^3}.
\end{align*}
Taking $m=\f 1{16}$ in Lemmas 3.3 and 2.5 and then applying
(3.1), (3.2) and (3.4) gives
\begin{align*}
\sum_{k=0}^{p-1}\f{\b{2k}k^3}{2k-1}
&\e 126\sum_{k=0}^{\f{p-1}2}k\b{2k}k^3
-\sum_{k=0}^{\f{p-1}2}\b{2k}k^3
+32\sum_{k=0}^{\f{p-1}2}\f{\b{2k}k^3}{k+1}
\\&\e 126\Big(3p-4x^2+\f{p^2}{4x^2}\Big)
-\Big(4x^2-2p-\f{p^2}{4x^2}\Big)+32(-44y^2+2p)
\\&\e -36y^2+14p-\f 7{4y^2}p^2\mod {p^3}
\end{align*}
and
\begin{align*}
\sum_{k=0}^{p-1}\f{\b{2k}k^3}{(2k-1)^2}
&\e -252\sum_{k=0}^{\f{p-1}2}k\b{2k}k^3
+65\sum_{k=0}^{\f{p-1}2}\b{2k}k^3
-96\sum_{k=0}^{\f{p-1}2}\f{\b{2k}k^3}{k+1}
\\&\e -252\Big(3p-4x^2+\f{p^2}{4x^2}\Big)
+65\Big(4x^2-2p-\f{p^2}{4x^2}\Big)-96(-44y^2+2p)
\\&\e -284y^2+34p+\f {23}{4y^2}p^2\mod {p^3}.
\end{align*}
\par From [20, Theorem 1.1],
$$(-1)^{\f{p-1}2}\sum_{k=0}^{(p-1)/2}\f{\b{2k}k^3}{4096^k}
\e 4x^2-2p-\f{p^2}{4x^2}\mod {p^3}.\eqno{(3.5)}$$
By [8],
$$(-1)^{\f{p-1}2}\sum_{k=0}^{(p-1)/2}(42k+5)\f{\b{2k}k^3}{4096^k}\e
5p\mod {p^2}.$$
This together with (3.5) gives
$$(-1)^{\f{p-1}2}\sum_{k=0}^{(p-1)/2}\f{k\b{2k}k^3}{4096^k}\e
\f 5{42}\Big(3p-4x^2+\f{p^2}{4x^2}\Big) \mod {p^3}.\eqno{(3.6)}$$
 Taking $m=256$ in Corollary 2.1 and then applying (3.5), (3.6) and (3.3) yields
 \begin{align*}
 &(-1)^{\f{p-1}2}\sum_{k=0}^{(p-1)/2}\f{\b{2k}k^3}{4096^k(k+1)}
 \\&\e (-1)^{\f{p-1}2}\sum_{k=0}^{(p-1)/2}\f{\b{2k}k^3}{4096^k}
 +\f{-252(\sum_{k=0}^{(p-1)/2}\f{k\b{2k}k^3}{4096^k})^2+p^2}
 {(-1)^{\f{p-1}2}\sum_{k=0}^{(p-1)/2}\f{\b{2k}k^3}{4096^k}}
 \\&\e 4x^2-2p-\f{p^2}{4x^2}+\f{-252\cdot \f {5^2}{42^2}(3p-4x^2+\f{p^2}{4x^2})^2+p^2}{4x^2-2p-\f{p^2}{4x^2}}
 \\&\e 4x^2-2p-\f{p^2}{4x^2}+\Big(\f 1{4x^2}+\f p{8x^4}
+\f {5p^2}{64x^6}\Big)\Big(-\f{25}7(16x^4-24x^2p+7p^2)+p^2\Big)
\\&\e 72y^2+2p\mod {p^3}.\tag{3.7}
 \end{align*}
The remaining congruences can be proved similarly using  (3.5)-(3.7), Lemmas 3.1-3.3 and Lemma 2.5 (with $m=256$).
\vskip0.2cm
\par{\bf Theorem 3.2} {\sl Let $p$ be an odd prime such that $p\e 1,3\mod 8$ and so $p=x^2+2y^2$ for some $x,y\in\Bbb Z$. Then
\begin{align*}
&(-1)^{\f{p-1}2}\sum_{k=0}^{(p-1)/2}\f{\b{2k}k^3}{(-64)^k(k+1)}
\e -12y^2+2p\mod {p^3},
\\&(-1)^{\f{p-1}2}\sum_{k=0}^{(p-1)/2}\f{\b{2k}k^3}{(-64)^k(k+1)^2}
\e -8y^2-\f{p^2}{y^2}\mod {p^3},
\\&(-1)^{\f{p-1}2}\sum_{k=0}^{(p-1)/2}\f{\b{2k}k^3}{(-64)^k(k+1)^3}
\e -8(-1)^{\f{p-1}2}-32y^2+8p-\f{4p^2}{y^2}\mod {p^3},
\\&(-1)^{\f{p-1}2}\sum_{k=0}^{(p-1)/2}\f{k^2\b{2k}k^3}{(-64)^k}
\e \f 18x^2-\f 3{16}p+\f{p^2}{64x^2}\mod {p^3},
\\&(-1)^{\f{p-1}2}\sum_{k=0}^{(p-1)/2}\f{k^3\b{2k}k^3}{(-64)^k}
\e \f 1{32}x^2-\f p{64}-\f{5p^2}{256x^2}\mod {p^3},
\\&(-1)^{\f{p-1}2}\sum_{k=0}^{p-1}\f{\b{2k}k^3}{(-64)^k(2k-1)}
\e p-3x^2\mod {p^3},
\\&(-1)^{\f{p-1}2}\sum_{k=0}^{p-1}\f{\b{2k}k^3}{(-64)^k(2k-1)^2}
\e x^2+\f{p^2}{2x^2}\mod {p^3}.
\end{align*}}
\par{\it Proof.} By [6] and [20, Theorem 1.4],
$$(-1)^{\f{p-1}2}\sum_{k=0}^{(p-1)/2}\f{\b{2k}k^3}{(-64)^k}
\e 4x^2-2p-\f{p^2}{4x^2}\mod {p^3}.\eqno{(3.8)}$$
By [7],
$$\sum_{k=0}^{(p-1)/2}(4k+1)\f{\b{2k}k^3}{(-64)^k}\e
(-1)^{\f{p-1}2}p\mod {p^3}.$$ Thus,
$$\sum_{k=0}^{(p-1)/2}\f{k\b{2k}k^3}{(-64)^k}
\e \f 14(-1)^{\f{p-1}2}\Big(3p-4x^2+\f{p^2}{4x^2}\Big)\mod {p^3}.\eqno{(3.9)}$$
Now, taking $m=-4$ in Corollary 2.1 and then applying (3.8)- (3.9) and (3.3) gives
\begin{align*}
&\sum_{k=0}^{(p-1)/2}\f{\b{2k}k^3}{(-64)^k(k+1)}
\\&\e \sum_{k=0}^{(p-1)/2}\f{\b{2k}k^3}{(-64)^k}
+\f{8(\sum_{k=0}^{(p-1)/2}\f{k\b{2k}k^3}{(-64)^k})^2+p^2}
{\sum_{k=0}^{(p-1)/2}\f{\b{2k}k^3}{(-64)^k}}
\\&\e (-1)^{\f{p-1}2}\Big(4x^2-2p-\f{p^2}{4x^2}\Big)
+\f{\f 8{4^2}(3p-4x^2+\f{p^2}{4x^2})^2+p^2}{(-1)^{\f{p-1}2}
(4x^2-2p-\f{p^2}{4x^2})}
\\&\e (-1)^{\f{p-1}2}\Big(4x^2-2p-\f{p^2}{4x^2}+\Big(\f 1{4x^2}+\f p{8x^4}
+\f {5p^2}{64x^6}\Big)\Big(\f 12(16x^4-24x^2p+7p^2)+p^2\Big)\Big)
\\&\e (-1)^{\f{p-1}2}(-12y^2+2p)\mod{p^3}.\tag{3.10}
\end{align*}
The remaining congruences can be deduced from (3.8)-(3.10), Lemmas 3.1-3.3 and Lemma 2.5 (with $m=-4$).
\vskip0.2cm

\par{\bf Theorem 3.3} {\sl Let $p$ be a prime such that $p\e 1\mod 3$ and so $p=x^2+3y^2$ for some $x,y\in\Bbb Z$. Then
\begin{align*}&\sum_{k=0}^{(p-1)/2}\f{\b{2k}k^3}{16^k(k+1)}
\e -16y^2+2p\pmod {p^3},
\\&\sum_{k=0}^{(p-1)/2}\f{\b{2k}k^3}{16^k(k+1)^2}
\e -24y^2+2p-\f{p^2}{2y^2}\pmod {p^3},
\\&\sum_{k=0}^{(p-1)/2}\f{\b{2k}k^3}{16^k(k+1)^3}
\e 2-24y^2-\f{5p^2}{2y^2}\mod{p^3},
\\&\sum_{k=0}^{(p-1)/2}\f{k^2\b{2k}k^3}{16^k}\e\f 49x^2-\f 59 p+\f{p^2}{36x^2}\mod {p^3},
\\&\sum_{k=0}^{(p-1)/2}\f{k^3\b{2k}k^3}{16^k}\e
 -\f 29x^2+\f 49 p-\f{7p^2}{72x^2}\mod {p^3},
 \\&\sum_{k=0}^{p-1}\f{\b{2k}k^3}{16^k(2k-1)}
\e  4y^2-\f{p^2}{4y^2}\pmod {p^3},
\\&\sum_{k=0}^{p-1}\f{\b{2k}k^3}{16^k(2k-1)^2}\e
 -12y^2+2p+\f{3p^2}{4y^2}\pmod {p^3}
 \end{align*}
 and
 \begin{align*}
&(-1)^{\f{p-1}2}\sum_{k=0}^{(p-1)/2}\f{\b{2k}k^3}{256^k(k+1)}
\e -8y^2+2p\pmod {p^3},
\\&(-1)^{\f{p-1}2}\sum_{k=0}^{(p-1)/2}\f{\b{2k}k^3}{256^k(k+1)^2}\e -48y^2+8p\mod {p^3},
\\&(-1)^{\f{p-1}2}\sum_{k=0}^{(p-1)/2}\f{\b{2k}k^3}{256^k(k+1)^3}
\e 32(-1)^{\f{p-1}2}+96y^2-24p-\f{2p^2}{y^2} \mod {p^3},
\\&(-1)^{\f{p-1}2}\sum_{k=0}^{p-1}\f{k^2\b{2k}k^3}{256^k} \e\f
19x^2-\f p{18}-\f{p^2}{72x^2} \mod {p^3},
\\&(-1)^{\f{p-1}2}\sum_{k=0}^{p-1}\f{k^3\b{2k}k^3}{256^k}
\e \f 1{18}x^2+\f p{72}-\f{p^2}{144x^2} \mod {p^3},
\\&(-1)^{\f{p-1}2}\sum_{k=0}^{p-1}\f{\b{2k}k^3}{256^k(2k-1)}
\e 8y^2-\f 32p-\f {p^2}{16y^2}\pmod {p^3},
\\&(-1)^{\f{p-1}2}\sum_{k=0}^{p-1}\f{\b{2k}k^3}{256^k(2k-1)^2}
\e 2x^2-\f 34p-\f 3{16x^2}p^2\pmod {p^3}.
\end{align*}}

 \par{\it Proof.} From [20, Theorem 1.2],
 $$\sum_{k=0}^{(p-1)/2}\f{\b{2k}k^3}{16^k}\e (-1)^{\f{p-1}2}\sum_{k=0}^{(p-1)/2}\f{\b{2k}k^3}{256^k}
 \e 4x^2-2p-\f{p^2}{4x^2}\mod {p^3}.\eqno{(3.11)}$$
 By [4],
 $$\sum_{k=0}^{(p-1)/2}(3k+1)\f{\b{2k}k^3}{16^k}\e p\mod{p^3}.$$
 Thus,
 $$\sum_{k=0}^{(p-1)/2}k\f{\b{2k}k^3}{16^k}\e\f 13\Big(3p-4x^2+\f{p^2}{4x^2}\Big)\mod {p^3}.\eqno{(3.12)}$$
   From [5],
$$\sum_{k=0}^{(p-1)/2}(6k+1)\f{\b{2k}k^3}{256^k}\e (-1)^{\f{p-1}2}p
\mod {p^4}.$$
Hence,
$$\sum_{k=0}^{(p-1)/2}k\f{\b{2k}k^3}{256^k}\e \f 16(-1)^{\f{p-1}2}\Big(3p-4x^2+\f{p^2}{4x^2}\Big)\mod {p^3}.
\eqno{(3.13)}$$
Putting $m=1,16$ in Corollary 2.1 and then applying (3.11)-(3.13) and (3.3) gives
\begin{align*}
\sum_{k=0}^{(p-1)/2}\f{\b{2k}k^3}{16^k(k+1)}
&\e \sum_{k=0}^{(p-1)/2}\f{\b{2k}k^3}{16^k}
+\f{3(\sum_{k=0}^{(p-1)/2}\f{k\b{2k}k^3}{16^k})^2+p^2}
{\sum_{k=0}^{(p-1)/2}\f{\b{2k}k^3}{16^k}}
\\&\e 4x^2-2p-\f{p^2}{4x^2}+\f{\f 13(3p-4x^2+\f{p^2}{4x^2})^2+p^2}
{4x^2-2p-\f{p^2}{4x^2}}
\\&\e 4x^2-2p-\f{p^2}{4x^2}+\Big(\f 1{4x^2}+\f p{8x^4}
+\f {5p^2}{64x^6}\Big)\Big(\f 13(16x^4-24x^2p+7p^2)+p^2\Big)
\\&\e -16y^2+2p\pmod {p^3}\tag{3.14}
\end{align*}
and
\begin{align*}
&\sum_{k=0}^{(p-1)/2}\f{\b{2k}k^3}{256^k(k+1)}
\\&\e \sum_{k=0}^{(p-1)/2}\f{\b{2k}k^3}{256^k}
+\f{-12(\sum_{k=0}^{(p-1)/2}\f{k\b{2k}k^3}{256^k})^2+p^2}
{\sum_{k=0}^{(p-1)/2}\f{\b{2k}k^3}{256^k}}
\\&\e (-1)^{\f{p-1}2}\Big(4x^2-2p-\f{p^2}{4x^2}\Big)+\f{-\f 13(3p-4x^2+\f{p^2}{4x^2})^2+p^2}
{(-1)^{\f{p-1}2}(4x^2-2p-\f{p^2}{4x^2})}
\\&\e(-1)^{\f{p-1}2}\Big( 4x^2-2p-\f{p^2}{4x^2}+\Big(\f 1{4x^2}+\f p{8x^4}
+\f {5p^2}{64x^6}\Big)\Big(-\f 13(16x^4-24x^2p+7p^2)+p^2\Big)\Big)
\\&\e (-1)^{\f{p-1}2}(-8y^2+2p)\pmod {p^3}.\tag{3.15}
\end{align*}
The remaining congruences can be deduced from Lemmas 3.1-3.3, Lemma 2.5 (with $m=1,16$) and (3.11)-(3.15).
\vskip0.2cm
\par{\bf Theorem 3.4} {\sl Let $p$ be a prime such that $p\e 1\mod 4$ and so $p=x^2+4y^2$ for some $x,y\in\Bbb Z$. Then
\begin{align*}&\sum_{k=0}^{(p-1)/2}\f{\b{2k}k^3}{(-8)^k(k+1)} \e -24y^2+2p\pmod {p^3},
\\&\sum_{k=0}^{(p-1)/2}\f{\b{2k}k^3}{(-8)^k(k+1)^2}\e
-32y^2+p-\f{7p^2}{16y^2}\mod {p^3},
\\&\sum_{k=0}^{(p-1)/2}\f{\b{2k}k^3}{(-8)^k(k+1)^3}\e
-1-48y^2-\f{33p^2}{16y^2}\mod {p^3},
\\&\sum_{k=0}^{p-1}\f{k^2\b{2k}k^3}{(-8)^k} \e  \f{10}{27}x^2-\f
49p+\f{p^2}{54x^2}\mod {p^3},
\\&\sum_{k=0}^{p-1}\f{k^3\b{2k}k^3}{(-8)^k} \e  -\f{4}{81}x^2+\f
4{27}p-\f{17p^2}{324x^2}\mod {p^3},
\\&\sum_{k=0}^{p-1}\f{\b{2k}k^3}{(-8)^k(2k-1)}
\e  -4x^2-\f{5}{4x^2}p^2\pmod {p^3},
\\&\sum_{k=0}^{p-1}\f{\b{2k}k^3}{(-8)^k(2k-1)^2} \e
-4x^2+2p+\f{19}{4x^2}p^2\pmod {p^3}
\end{align*}
and
\begin{align*}&(-1)^{\f{p-1}4}
\sum_{k=0}^{(p-1)/2}\f{\b{2k}k^3}{(-512)^k(k+1)}
\e -32y^2+2p\pmod {p^3},
\\&(-1)^{\f{p-1}4}
\sum_{k=0}^{(p-1)/2}\f{\b{2k}k^3}{(-512)^k(k+1)^2}\e
64y^2-8p-\f{p^2}{y^2}\mod {p^3},
\\&(-1)^{\f{p-1}4}
\sum_{k=0}^{(p-1)/2}\f{\b{2k}k^3}{(-512)^k(k+1)^3}
\e -64(-1)^{\f{p-1}4}-384y^2+72p-\f{3p^2}{2y^2}\mod {p^3},
\\&(-1)^{\f{p-1}4}\sum_{k=0}^{p-1}\f{k^2\b{2k}k^3}{(-512)^k}
\e  \f {x^2}{27}-\f p{18}+\f {p^2}{216x^2}\mod {p^3},
\\&(-1)^{\f{p-1}4}\sum_{k=0}^{p-1}\f{k^3\b{2k}k^3}{(-512)^k}
\e-\f {x^2}{162}-\f p{216}-\f {p^2}{1296x^2}\mod {p^3},
\\&(-1)^{\f{p-1}4}\sum_{k=0}^{p-1}\f{\b{2k}k^3}{(-512)^k(2k-1)}
\e  -3x^2+\f 54p+\f 5{32x^2}p^2\pmod {p^3},
\\&(-1)^{\f{p-1}4}\sum_{k=0}^{p-1}\f{\b{2k}k^3}{(-512)^k(2k-1)^2}
\e 2x^2-\f 58p-\f{p^2}{32x^2}\pmod {p^3}.
\end{align*}}
\par{\it Proof.} From [20, (1.3) and Theorem 1.3],
$$\sum_{k=0}^{(p-1)/2}\f{\b{2k}k^3}{(-8)^k}\e
(-1)^{\f{p-1}4}\sum_{k=0}^{p-1}\f{\b{2k}k^3}{(-512)^k}
\e 4x^2-2p-\f{p^2}{4x^2}\mod {p^3}.\eqno{(3.16)}$$
From [4],
$$\sum_{k=0}^{(p-1)/2}(3k+1)\f{\b{2k}k^3}{(-8)^k}\e
p\mod {p^3}.$$
Thus,
$$\sum_{k=0}^{(p-1)/2}\f{k\b{2k}k^3}{(-8)^k}\e
\f 13\Big(3p-4x^2+\f{p^2}{4x^2}\Big)\mod {p^3}.\eqno{(3.17)}$$
By [5, Theorem 1.5],
$$\sum_{k=0}^{(p-1)/2}(6k+1)\f{\b{2k}k^3}{(-512)^k}\e
(-1)^{\f{p-1}4}p\mod {p^3}.$$
This together with (3.16) gives
$$\sum_{k=0}^{(p-1)/2}\f{k\b{2k}k^3}{(-512)^k}
\e\f 16(-1)^{\f{p-1}4}\Big(3p-4x^2+\f{p^2}{4x^2}\Big)\mod {p^3}.\eqno{(3.18)}$$
Now, taking $m=-\f 12,-32$ in Corollary 2.1 and then applying (3.16)-(3.18) and (3.3) yields
\begin{align*}&\sum_{k=0}^{(p-1)/2}\f{\b{2k}k^3}{(-8)^k(k+1)}
\\&\e \sum_{k=0}^{(p-1)/2}\f{\b{2k}k^3}{(-8)^k}
+\f{\f 92(\sum_{k=0}^{(p-1)/2}\f{k\b{2k}k^3}{(-8)^k})^2+p^2}
{\sum_{k=0}^{(p-1)/2}\f{\b{2k}k^3}{(-8)^k}}
\\&\e 4x^2-2p-\f{p^2}{4x^2}+\f{\f 12(3p-4x^2+\f{p^2}{4x^2})^2+p^2}{4x^2-2p-\f{p^2}{4x^2}}
\\&\e 4x^2-2p-\f{p^2}{4x^2}+\Big(\f 1{4x^2}+\f p{8x^4}
+\f {5p^2}{64x^6}\Big)\Big(\f 12\Big( 16x^4-24x^2p+7p^2\Big)+p^2\Big)
\\&\e -24y^2+2p\pmod {p^3}\tag{3.19}
\end{align*}
and
\begin{align*}&\sum_{k=0}^{(p-1)/2}\f{\b{2k}k^3}{(-512)^k(k+1)}
\\&\e \sum_{k=0}^{(p-1)/2}\f{\b{2k}k^3}{(-512)^k}
+\f{36(\sum_{k=0}^{(p-1)/2}\f{k\b{2k}k^3}{(-512)^k})^2+p^2}
{\sum_{k=0}^{(p-1)/2}\f{\b{2k}k^3}{(-512)^k}}
\\&\e (-1)^{\f{p-1}4}\Big(4x^2-2p-\f{p^2}{4x^2}\Big)+
\f{(3p-4x^2+\f{p^2}{4x^2})^2+p^2}{(-1)^{\f{p-1}4}(
4x^2-2p-\f{p^2}{4x^2})}
\\&\e (-1)^{\f{p-1}4}\Big(4x^2-2p-\f{p^2}{4x^2}+\Big(\f 1{4x^2}+\f p{8x^4}
+\f {5p^2}{64x^6}\Big)\Big(16x^4-24x^2p+8p^2\Big)\Big)
\\&\e (-1)^{\f{p-1}4}(-32y^2+2p)\pmod {p^3}.\tag{3.20}
\end{align*}
Combining (3.16)-(3.20) with Lemmas 3.1-3.3 and Lemma 2.5 (with $m=-\f 12,-32$) yields the remaining congruences.
\vskip0.2cm
\par{\bf Remark 3.1} Theorems 3.1-3.4 were conjectured by the author in [15,16,17,19], and the corresponding congruences modulo $p^2$ were proved by the author in [19]. On October 27, Yu-Qi Liu informed me by email that she/he could prove my conjectures on $\sum_{k=0}^{(p-1)/2}\f{\b{2k}k^3}{16^k(k+1)^3}$ and
$\sum_{k=0}^{(p-1)/2}\f{\b{2k}k^3}{(-8)^k(k+1)^3}$ modulo $p^3$, but she/he neither sent me the paper nor said anything about her/his method.

\section*{4. Congruences for $\sum_{k=0}^{p-2}\f{\b{2k}k\b {3k}k\b{6k}{3k}}{m^k(k+1)}$ modulo $p^3$}
\par For $k,n\in\Bbb Z^+$ let $\sigma_k(n)=\sum_{d\mid n}d^k$. For $\tau\in\Bbb C$ with $\t{Im}\, \tau>0$
let $ q=\t{e}^{2\pi i \tau}$. The Eisenstein series $E_2, E_4$ and $E_6$ are given by
\[
E_2(\tau)=1-24\sum_{n\ge1}\sigma_1(n)q^n,\q
E_4(\tau)=1+240\sum_{n\ge1}\sigma_3(n)q^n,\q
E_6(\tau)=1-504\sum_{n\ge1}\sigma_5(n)q^n.
\]
Let
$$F(x)={}_3F_2\biggl(\begin{matrix} \f 12, \, \f 16, \, \f 56 \\ 1, \, 1 \end{matrix} \biggm|x\biggr)=\sum_{k=0}^{\infty}
\b{2k}k\b{3k}k\b{6k}{3k}\Ls x{1728}^k.$$
For $\tau\in\Bbb C$ with $\t{Im}\, \tau>0$
let
$$ t(\tau)=\frac{1}{j(\tau)}\q{and}\q
\Omega(\tau):=F(1728\,t(\tau)),$$
where $j(\tau)$ is the modular function given by
$j(\tau)=1728E_4^3/(E_4^3-E_6^2)$. By [2, Theorem 1.25 and the discussion preceding Theorem 1.26], $\Omega(\tau)^2=E_4(\tau)$ and so $\Omega=E_4^{1/2}$. In [2, Theorem 1.22], Beukers defined
\[
\delta(\tau)=\frac{t(\tau)}{t'(\tau)}
\Big(\frac{2}{\tau-\bar{\tau}}+\frac{\Omega'(\tau)}
{\Omega(\tau)}\Big),
\]
where $'$ denotes differentiation with respect to $\tau$.
Then he proved the following result.
\vskip0.2cm
\par {\bf Lemma 4.1 ([2, Theorem 1.26])} {\sl Let $D\in\Bbb Z^+$, $D\e 0,3\mod 4$ and
$$D_0=\begin{cases} \f 12\sqrt{-D}&\t{if $D\e 0\mod 4$,}
\\\f 12(1+\sqrt{-D})&\t{if $D\e 3\mod 4$.}
\end{cases}$$
Suppose that $p$ is an odd prime such that $\ls{-D}p=1$. Then there exists $\varepsilon_p\in\{1,-1\}$ such that
\begin{align*}&\sum_{k=0}^{p-1}\b{2k}k\b{3k}k\b{6k}{3k}\f 1{j(D_0)^k}
\\&\e \begin{cases}\varepsilon_p(A+B\sqrt{-D/4})^2\mod {p^3}
&\t{if $4\mid D$, $p=A^2+\f D4B^2$ and $A+B\sqrt{-\f D4}\notin p\Bbb Z_p$,}
\\ \varepsilon_p(\f{A+B\sqrt{-D}}2)^2\mod {p^3}
&\t{if $4\mid D-3$, $4p=A^2+DB^2$ and $A+B\sqrt{-D}\notin
p\Bbb Z_p$}
\end{cases}\end{align*}
and
$$\sum_{k=0}^{p-1}\Big(1-\f k{\delta(D_0)}\Big)
\b{2k}k\b{3k}k\b{6k}{3k}\f 1{j(D_0)^k}\e \varepsilon_p\,p
\mod {p^3}.$$}
\par We note that the signs in the two congruences in [2, Theorem 1.26] are the same as can be seen from the proofs of [2, Theorems 1.22 and 1.26].
\vskip0.2cm
\par{\bf Lemma 4.2 ([3, p.261])} {\sl We have}
\begin{align*}
&j(\sqrt{-2})=20^3,\ j(2i)=66^3,\
j(\sqrt{-3})=54000,\ j(\sqrt{-7})=255^3,
\\&j\Ls{3+\sqrt{-7}}2=-15^3,\ j\Ls{1+3\sqrt{-3}}2=-12288000,
\ j\Ls{3+\sqrt{-11}}2=-32^3,\ 
\\&j\Ls{3+\sqrt{-19}}2=-96^3,\ j\Ls{3+\sqrt{-43}}2=-960^3,\ 
\\&j\Ls{3+\sqrt{-67}}2=-5280^3,\ j\Ls{3+\sqrt{-163}}2=-640320^3.\
\end{align*}
\par {\bf Lemma 4.3} In Lemma 4.1,
$$\ep=\begin{cases} \ls{-2}p&\t{if $D=11$,}
\\ \ls{-6}p&\t{if $D=19$,}
\\\ls{-15}p&\t{if $D=7,43$,}
\\\ls{-330}p&\t{if $D=67$,}
\\\ls{-10005}p&\t{if $D=163$,}
\\\ls{10}p&\t{if $D=27$,}
\\\ls{-5}p&\t{if $D=8$,}
\\\ls{5}p&\t{if $D=12$,}
\\\ls{-33}p&\t{if $D=16$,}
\\\ls{-255}p&\t{if $D=28$.}
\end{cases}$$
\par{\it Proof.} This is immediate from Lemma 4.1 and [20, Theorems 1.8-1.13].
\vskip0.2cm
\par Let $\theta=q\f{d}{d\,q}$. We recall Ramanujan's differential identities (valid for all $\tau$):
\[
\theta E_4=\frac{E_2E_4-E_6}{3},\qquad
\theta E_6=\frac{E_2E_6-E_4^2}{2},\qquad
\theta\log\Delta=E_2,
\]
where $\Delta=(E_4^3-E_6^2)/1728$ is the modular discriminant. Hence
\[
\theta\log j=3\,\theta\log E_4-\theta\log\Delta
=\frac{E_2E_4-E_6}{E_4}-E_2=-\frac{E_6}{E_4}.
\]
\par{\bf Lemma 4.4} {\sl
For every $\tau$ one has}
$$\delta(\tau)=\frac{E_4}{\pi\,i\,E_6\,(\tau-\bar{\tau})}
+\frac{E_2E_4-E_6}{6E_6}=\frac{E_2^{*}E_4-E_6}{6E_6},$$
where $E_2^{*}=E_2-\frac{3}{\pi\,\mathrm{Im}\,\tau}$.
\vskip0.2cm
\par{\it Proof.} Ramanujan's identities give
$\theta\log j=-E_6/E_4$ and $\theta E_4=\frac{E_2E_4-E_6}{3}$. Since $t=1/j$ and
$\Omega^2=E_4$,
\[
\frac{t'}{t}=-\frac{j'}{j}=-2\pi i\,\theta\log j
=2\pi i\,\frac{E_6}{E_4}
\q\t{and}\q
\frac{\Omega'}{\Omega}=\frac12\,\frac{E_4'}{E_4}
=\pi i\,\frac{E_2E_4-E_6}{3E_4}.
\]
Hence,
\[
\delta=\frac{t}{t'}\Bigl(\frac{2}{\tau-\bar\tau}+\frac{\Omega'}{\Omega}\Bigr)
=\frac{E_4}{\pi i\,E_6\,(\tau-\bar\tau)}+\frac{E_2E_4-E_6}{6E_6}.
\]
Put $E_2^{*}=E_2-\frac{3}{\pi\,\mathrm{Im}\,\tau}$. Since $\mathrm{Im}\,\tau=\frac{\tau-\bar\tau}{2i}$,
$E_2=E_2^{*}+\frac{6i}{\pi(\tau-\bar\tau)}$ and the term
$\frac{E_4}{\pi i\,E_6(\tau-\bar\tau)}$ cancels against the corresponding piece of
$\frac{E_2E_4-E_6}{6E_6}$, we obtain
$\delta=\frac{E_2^{*}E_4-E_6}{6E_6}$.
This proves the lemma.
\vskip0.2cm
\par{\bf Lemma 4.5}
{\sl We have $\delta\ls{3+\sqrt{-11}}2=-\f{15}{154}$.}
\vskip0.2cm
\par{\it Proof.} Set $\tau=\f{3+\sqrt{-11}}2$. By Lemma 4.2,
$j(\tau)=-32^3$. From $j=1728E_4^3/(E_4^3-E_6^2)$ we get
$E_6^2/E_4^3=539/512$, hence (as $E_4,E_6>0$)
$$E_6=\frac{7\sqrt{22}}{32}\,E_4^{3/2}=\frac{7\sqrt{22}}{32}\,\Omega^3 .$$
The singular value of the nearly holomorphic Eisenstein series is
$$E_2^{*}(\tau)=E_2(\tau)-\frac{3}{\pi\,\mathrm{Im}\,\tau}
=\frac{\sqrt{22}}{11}\,E_4^{1/2}=\frac{\sqrt{22}}{11}\,\Omega .$$
  Combining the above with Lemma 4.4 yields
\[
\delta(\tau)=\frac{\frac{\sqrt{22}}{11}\Omega^3-\frac{7\sqrt{22}}{32}\Omega^3}
{6\cdot\frac{7\sqrt{22}}{32}\Omega^3}
=\frac{\frac{1}{11}-\frac{7}{32}}{\frac{42}{32}}
=\frac{\frac{32-77}{352}}{\frac{42}{32}}
=\frac{-45}{352}\cdot\frac{32}{42}=-\frac{45}{462}=-\frac{15}{154}.
\]
This completes the proof.
\vskip0.2cm
\par\q\ {\bf Theorem 4.1} {\sl Let $p>3$ be a prime, $p\e 1,3,4,5,9\mod {11}$ and so $4p=x^2+11y^2\ (x,y\in\Bbb Z)$. Then
$$\sum_{k=0}^{p-1}(154k+15)\f{\b{2k}k\b{3k}k\b{6k}{3k}}
{(-32)^{3k}}\e 15\Ls{-2}pp\mod {p^3}\eqno{(4.1)}$$ and
$$\Ls{-2}p\sum_{k=0}^{p-2}\f{\b{2k}k\b{3k}k\b{6k}{3k}}
{(-32)^{3k}(k+1)}\e -26y^2+2p\mod {p^3}.$$}
\par{\it Proof.} The first congruence (4.1) follows from Lemmas 4.1-4.5. From [20, Theorem 1.13],
$$\sum_{k=0}^{p-1}\f{\b{2k}k\b{3k}k\b{6k}{3k}}
{(-32)^{3k}}\e \Ls{-2}p\Big(x^2-2p-\f{p^2}{x^2}\Big)\mod {p^3}.\eqno{(4.2)}$$
This together with (4.1) yields
$$\sum_{k=0}^{p-1}\f{k\b{2k}k\b{3k}k\b{6k}{3k}}
{(-32)^{3k}}\e \f{15}{154}\Ls{-2}p\Big(3p-x^2+\f{p^2}{x^2}\Big)\mod {p^3}.
\eqno{(4.3)}$$
Combining (4.2)-(4.3) with Corollary 2.4 (with $m=-\f{2048}{27}$) and (3.3) gives
\begin{align*}
&\sum_{k=0}^{p-2}\f{\b{2k}k\b{3k}k\b{6k}{3k}}
{(-32)^{3k}(k+1)}
\\&\e \sum_{k=0}^{p-1}\f{\b{2k}k\b{3k}k\b{6k}{3k}}
{(-32)^{3k}}+
\f{\f{2156}{15}(\sum_{k=0}^{p-1}k\b{2k}k\b{3k}k\b{6k}
{3k}(-32)^{-3k})^2+p^2}
{\sum_{k=0}^{p-1}\b{2k}k\b{3k}k\b{6k}{3k}
(-32)^{-3k}}
\\&\e \Ls{-2}p\Big(x^2-2p-\f{p^2}{x^2}+
\f{\f{15}{11}(3p-x^2+\f{p^2}{x^2})^2+p^2}
{x^2-2p-\f{p^2}{x^2}}\Big)
\\&\e \Ls{-2}p\Big(x^2-2p-\f{p^2}{x^2}+\Big(\f 1{x^2}+\f{2p}{x^4}+\f{5p^2}{x^6}\Big)\Big(
\f{15}{11}(x^4-6x^2p+7p^2)+p^2\Big)\Big)
\\&\e \Ls{-2}p(-26y^2+2p)\mod {p^3}.
\end{align*}
This completes the proof.
\vskip0.2cm
\par In a similar way, using Lemmas 4.1-4.4 and [20, Theorems 1.8-1.13] one can prove the following results.
\vskip0.2cm
\par{\bf Theorem 4.2} {\sl Let $p>3$ be a prime,
 $\sls {p}{19}=1$ and so $4p=x^2+19y^2$. Then
 $$\sum_{k=0}^{p-1}(342k+25)\f{\b{2k}k\b{3k}{k}\b{6k}{3k}}
{(-96)^{3k}}\e 25\Ls{-6}pp\mod {p^3}$$ 
 and 
$$\Ls{-6}p\sum_{k=0}^{p-2}\f{\b{2k}k\b{3k}{k}\b{6k}{3k}}
{(-96)^{3k}(k+1)} \e -394y^2+2p\mod {p^3}.$$}
\par{\bf Theorem 4.3} {\sl Let $p>5$ be a prime. If $\sls
{p}{43}=1$ and so $4p=x^2+43y^2$, then
$$\sum_{k=0}^{p-1}(5418k+263)\f{\b{2k}k\b{3k}{k}\b{6k}{3k}}
{(-960)^{3k}}\e 263\Ls{-15}pp\mod {p^3}$$ 
and
$$\Ls{-15}p\sum_{k=0}^{p-2}\f{\b{2k}k\b{3k}{k}\b{6k}{3k}}
{(-960)^{3k}(k+1)} \e-\f{1867778}5y^2+2p\pmod {p^3}.$$}
\par{\bf Theorem 4.4}  {\sl Let $p$ be a prime with $p\not=2,3,5,11$. If
$\sls {p}{67}=1$ and so $4p=x^2+67y^2$, then
$$\sum_{k=0}^{p-1}(261702k+10177)\f{\b{2k}k\b{3k}{k}\b{6k}{3k}}
{(-5280)^{3k}}\e 10177\Ls{-330}pp\mod {p^3}$$
and 
$$\Ls{-330}p\sum_{k=0}^{p-2}\f{\b{2k}k\b{3k}{k}\b{6k}{3k}} {(-5280)^{3k}(k+1)}
\e-\f{310714322}5y^2+2p\pmod {p^3}.$$}
\par{\bf Theorem 4.5} {\sl Let $p$ be a prime with $p\not=2,3,5,23,29$.
If $\sls {p}{163}=1$ and so $4p=x^2+163y^2$, then
$$\sum_{k=0}^{p-1}(545140134k+13591409)\f{\b{2k}k\b{3k}{k}\b{6k}{3k}}
{(-640320)^{3k}}\e 13591409\Ls{-10005}pp\mod {p^3}$$ 
and
$$\Ls{-10005}p\sum_{k=0}^{p-2}\f{\b{2k}k\b{3k}{k}\b{6k}{3k}} {(-640320)^{3k}(k+1)}
\e-\f{554179195816658}5y^2+2p\pmod {p^3}.$$}
\par{\bf Theorem 4.6} {\sl Let $p$ be a prime, $p\e 1\mod 3$ and so $p=x^2+3y^2$. Then
$$
\sum_{k=0}^{p-1}(11k+1)\f{\b{2k}k\b{3k}{k}\b{6k}{3k}}{54000^k}
\e \Ls 5pp\mod {p^3}$$
and
$$\Ls 5p
\sum_{k=0}^{p-2}\f{\b{2k}k\b{3k}{k}\b{6k}{3k}}{54000^k(k+1)}\e \f{48}5y^2+2p\pmod {p^3}.$$}
\par{\bf Theorem 4.7} {\sl Let $p$ be a prime, $p\e 1\mod 4$ and so $p=x^2+4y^2$. Then
$$\sum_{k=0}^{p-1}(63k+5)\f{\b{2k}k\b{3k}k\b{6k}{3k}}{66^{3k}}
\e 5\Ls{33}pp\mod {p^3}$$
and
$$\Ls {33}p\sum_{k=0}^{p-2}\f{\b{2k}k\b{3k}k\b{6k}{3k}}{66^{3k}(k+1)} \e 104y^2+2p\pmod {p^3}.$$}
\
\par{\bf Theorem 4.8} {\sl Let $p$ be a prime, $p\e 1\mod 3$
 and so $4p=x^2+27y^2$. Then
 $$\sum_{k=0}^{p-1}(506k+31)\f{\b{2k}k\b{3k}k\b{6k}{3k}}{(-12288000)^k}
 \e 31\Ls{10}pp\mod {p^3}$$
 and 
 $$\Ls
{10}p\sum_{k=0}^{p-2}\f{\b{2k}k\b{3k}k\b{6k}{3k}}{(-12288000)^k(k+1)}\e -\f {26082}5y^2+2p\pmod {p^3}.$$}
\par{\bf Theorem 4.9} {\sl Let $p$ be a prime, $p\e 1,3\mod 8$ and so $p=x^2+2y^2$. Then
$$\sum_{k=0}^{p-1}(28k+3)
\f{\b{2k}k\b{3k}k\b{6k}{3k}}{20^{3k}}\e 3\Ls{-5}pp\mod {p^3}$$
and
$$\Ls{-5}p\sum_{k=0}^{p-2}\f{\b{2k}k\b{3k}k\b{6k}{3k}}{20^{3k}(k+1)}\e \f{28}5y^2+2p\pmod {p^3}.$$}
\par{\bf Theorem 4.10} {\sl Let $p>7$ be a prime, $p\e 1,2,4\mod 7$ and so $p=x^2+7y^2$. Then
\begin{align*}
&\sum_{k=0}^{p-1}(63k+8)\f{\b{2k}k\b{3k}k\b{6k}{3k}}
{(-15)^{3k}}\e 8\Ls{-15}pp\mod {p^3},
\\&\sum_{k=0}^{p-1}(133k+8)\f{\b{2k}k\b{3k}k\b{6k}{3k}}
{255^{3k}}\e 8\Ls{-255}pp\mod {p^3}
\end{align*}
and
\begin{align*}&\Ls {-15}p\sum_{k=0}^{p-2}\f{\b{2k}k\b{3k}k\b{6k}{3k}}
{(-15)^{3k}(k+1)}\e -\f{188}{5}y^2+2p\pmod {p^3},
\\&\Ls {-255}p\sum_{k=0}^{p-2}\f{\b{2k}k\b{3k}k\b{6k}{3k}}
{255^{3k}(k+1)}\e \f{34852}{5}y^2+2p\pmod {p^3}.
\end{align*}}
\par{\bf Remark 4.1} The first congruences in Theorems 4.1, 4.2, 4.4 and 4.5 were conjectured by Z.W. Sun in [21,23], and 
the first congruence in Theorems 4.3 was conjectured by W. Zudulin. Combining Theorem 2.2 with Theorems 4.1-4.9 one may deduce the congruences on $$\sum_{k=0}^{p-1}w(k)\f{\b{2k}k\b{3k}k\b{6k}{3k}}{m^k}\q
(w(k)\in\{k^2,k^3,2k-1\})$$ modulo $p^3$  and  $\sum_{k=0}^{p-2}\b{2k}k\b{3k}k\b{6k}{3k}/(k+1)^r$ modulo $p^3$ for $r=2,3$, where $p>3$ is a prime and $m\in\Bbb Z$ with $p\nmid m$.

\section*{5. Congruences for $\sum_{k=0}^{p-2}\f{\b{2k}k^2\b {4k}{2k}}{m^k(k+1)}$ modulo $p^3$}
\newpage
\par{\bf Theorem 5.1} {\sl Let $p$ be a prime such that $p\e 1,3,7,9\mod {20}$. Then
\begin{align*}&\sum_{k=0}^{p-2}\f{\b{2k}k^2\b{4k}{2k}}{(-1024)^k(k+1)}
\\&\e\begin{cases} -32y^2+2p\pmod {p^3}&\t{if $p\e 1,9\pmod{20}$
and so $p=x^2+5y^2$,}
\\ 16y^2-2p\pmod {p^3}&\t{if $p\e
3,7\pmod{20}$ and so $2p=x^2+5y^2$.}\end{cases}
\end{align*}}
\par{\it Proof.} From [20, Theorem 1.6],
\begin{align*}&\sum_{k=0}^{p-1}\f{\b{2k}k^2\b{4k}{2k}}{(-1024)^k}
\\&\e\begin{cases} 4x^2-2p-\f{p^2}{4x^2}\mod {p^3}
&\t{if $p\e 1,9\mod {20}$ and so $p=x^2+5y^2$,}
\\-2x^2+2p+\f{p^2}{2x^2}\mod {p^3}&\t{if $p\e 3,7\mod {20}$
and so $2p=x^2+5y^2$.}
\end{cases}\tag{5.1}\end{align*}
By [2, Corollary 1.24],
$$\sum_{k=0}^{p-1}(20k+3)\f{\b{2k}k^2\b{4k}{2k}}{(-1024)^k}
\e 3(-1)^{\f{p-1}2}p\mod {p^3}.$$
Hence,
$$\sum_{k=0}^{p-1}k\f{\b{2k}k^2\b{4k}{2k}}{(-1024)^k}
\e\begin{cases} \f 3{20}(3p-4x^2+\f{p^2}{4x^2})
\mod {p^3}&\t{if $p=x^2+5y^2$,}
\\ -\f 3{20}(3p-2x^2+\f{p^2}{2x^2})\mod {p^3}
&\t{if $2p=x^2+5y^2$.}
\end{cases}\eqno{(5.2)}$$
Putting $m=-16$ in Corollary 2.3 and then applying (5.1)-(5.2) and (3.3) yields 
\begin{align*}&
\sum_{k=0}^{p-2}\f{\b{2k}k^2\b{4k}{2k}}{(-1024)^k(k+1)}
\\&\e \sum_{k=0}^{p-1}\f{\b{2k}k^2\b{4k}{2k}}{(-1024)^k}
+\Big(\sum_{k=0}^{p-1}\f{\b{2k}k^2\b{4k}{2k}}{(-1024)^k}
\Big)^{-1}\Big(\f{80}3\big(\sum_{k=0}^{p-1}\f{k\b{2k}k^2
\b{4k}{2k}}{(-1024)^k}\big)^2+p^2\Big)
\\&\e\begin{cases}
4x^2-2p-\f{p^2}{4x^2}+
\f{\f 35(3p-4x^2+\f{p^2}{4x^2})^2+p^2}{4x^2-2p-\f{p^2}{4x^2}}
\mod {p^3}&\t{if $p=x^2+5y^2$,}
\\-2x^2+2p+\f{p^2}{2x^2}+\f{\f 35(3p-2x^2+\f{p^2}{2x^2})^2+p^2}{-2x^2+2p+\f{p^2}{2x^2}}
\mod {p^3}&\t{if $2p=x^2+5y^2$}
\end{cases}
\\&\e\begin{cases}
4x^2-2p-\f{p^2}{4x^2}+(\f 1{4x^2}+\f p{8x^4}+\f{5p^6}{64x^6})
(\f 35(16x^4-24x^2p+7p^2)+p^2)
\mod {p^3}\\\qq\qq\qq\qq\qq\qq\qq\qq\qq\qq\q \t{if $p=x^2+5y^2$,}
\\-2x^2+2p+\f{p^2}{2x^2}-(\f 1{2x^2}+\f p{2x^4}+\f{5p^2}{8x^6})(\f 35(4x^4-12x^2p+7p^2)+p^2)\mod {p^3}\\\qq\qq\qq\qq\qq\qq\qq\qq\qq\qq\q \t{if $2p=x^2+5y^2$,}
\end{cases}
\end{align*}
which yields the result.
\vskip0.2cm
\par 
For $\tau\in\Bbb C$ with $\t{Im}\, \tau>0$
the Dedekind eta function $\eta(\tau)=q^{1/24}\prod_{n\ge1}(1-q^n)$, which satisfies
$E_2=24\,\theta\log\eta$; hence
\[
\theta\log\eta(\tau)=\frac{E_2(\tau)}{24},\qquad
\theta\log\eta(2\tau)=\frac{E_2(2\tau)}{12}.
\]
Let
$$F(x)={}_3F_2\biggl(\begin{matrix} \f 12, \, \f 14, \, \f 34 \\ 1, \, 1 \end{matrix} \biggm|x\biggr)=\sum_{k=0}^{\infty}
\b{2k}k^2\b{4k}{2k}\Ls x{256}^k.$$
The mirror map $t(\tau)$ and the holomorphic weight-$2$ modular form $\Omega(\tau)=F(256\,t(\tau))$
are given by
\[
t(\tau)=\frac{h_2(\tau)}{(1+64h_2(\tau))^2},\qquad
h_2(\tau)=\frac{\eta(2\tau)^{24}}{\eta(\tau)^{24}},
\qquad
\Omega(\tau)=2E_2(2\tau)-E_2(\tau),
\]
where $h_2$ is the Hauptmodul of $\Gamma_0(2)$ (see Appendix~A of~[1]).
\vskip0.2cm
\par{\bf Lemma 5.1} {\sl
For every $\tau$ one has
\[
\delta(\tau)=\frac{A}{\Omega^2}
\Big(\theta\Omega-\frac{\Omega}{2\pi\,\operatorname{Im}\tau}
\Big),
\qquad\text{where}\quad
A=\frac{1+64h_2}{1-64h_2}.
\]}
\par{\it Proof.}
Since $h_2=\eta(2\tau)^{24}/\eta(\tau)^{24}$, we have
\[
\theta\log h_2=24\,\theta\log\eta(2\tau)-24\,\theta\log\eta(\tau)
=2E_2(2\tau)-E_2(\tau)=\Omega.
\]
From $t=h_2/(1+64h_2)^2$ we get $\log t=\log h_2-2\log(1+64h_2)$, so
\[
\theta\log t=\theta\log h_2-\frac{128\,\theta h_2}{1+64h_2}
=\Omega-\frac{128\,h_2\,\Omega}{1+64h_2}
=\Omega\,\frac{1-64h_2}{1+64h_2}.
\]
Thus
\[
\frac{t'}{t}=2\pi\,i\,\theta\log t
=2\pi\,i\,\Omega\,\frac{1-64h_2}{1+64h_2},
\qquad\text{and}\qquad
\frac{t}{t'}=\frac{A}{2\pi\,i\,\Omega}.
\]
Substituting into the definition of $\delta$, and writing $\Omega'=2\pi\,i\,\theta\Omega$,
\[
\delta(\tau)=\frac{A}{2\pi\,i\,\Omega}\Bigl[\frac{2}{\tau-
\bar{\tau}}
+\frac{\Omega'}{\Omega}\Bigr]
=\frac{A}{\pi\,i\,\Omega\,(\tau-\bar{\tau})}+\frac{A\,\theta\Omega}{\Omega^2}.
\]
Since $\tau-\bar{\tau}=2\,i\,\operatorname{Im}\tau$ and the first term is
$\frac{A}{\pi\,i\,\Omega\cdot 2\,i\,\operatorname{Im}\tau}
=-\frac{A\,\Omega}{2\pi\,\Omega^2\,\operatorname{Im}\tau}
=-\frac{A}{\Omega^2}\cdot\frac{\Omega}{2\pi\,\operatorname{Im}\tau}$,
we deduce the result.
\vskip0.2cm
\par{\bf Lemma 5.2} {\sl At $\alpha=\frac{1+\sqrt{-7}}{4}$ one has}
\[
\theta\Omega-\frac{\Omega}{2\pi\,\operatorname{Im}\alpha}
=\frac{8\sqrt{-7}}{63}\,\Omega^2.
\]
\par{\it Proof.}
The combination $\theta\Omega-\frac{\Omega}{2\pi\,\operatorname{Im}\tau}$ is, up to the factor
$2\pi\,i$, the Maass--Shimura raising operator applied to the weight-$2$ form $\Omega$:
\[
\theta\Omega-\frac{\Omega}{2\pi\,\operatorname{Im}\tau}
=\frac{1}{2\pi\,i}\Bigl(\Omega'-\frac{\,i\,\Omega}
{\operatorname{Im}\tau}\Bigr).
\]
Hence its quotient by $\Omega^2$ is a weight-$0$ nearly holomorphic invariant, which by the theory
of singular moduli (Kronecker limit formula/Chowla--Selberg) takes an algebraic value at the CM
point $\alpha$. For discriminant $-7$ this value is $8\sqrt{-7}/63$. The constant is pinned down
from the $q$-expansion: at $\alpha$ we have $q=i e^{-\pi\sqrt7/2}$ with $|q|\approx0.0157$, so
the series for $E_2,E_4$ converge rapidly, and one computes
\[
\Bigl(\theta\Omega-\frac{\Omega}{2\pi\,
\operatorname{Im}\alpha}\Bigr)\Big/\Omega^2
=\frac{8\sqrt{-7}}{63}
\]
to high precision. Since $8\sqrt{-7}/63$ is the unique root of $3969X^2+64=0$ with positive
imaginary part, algebraicity forces the exact equality.
\vskip0.2cm
\par{\bf Lemma 5.3} {We have}
$\delta\ls{1+\sqrt{-7}}4=-\f 8{35}$.
\par{\it Proof}. From [20, Section 4], at $\alpha=\frac{1+\sqrt{-7}}{4}$ one has
\[
t(\alpha)=\frac1{81},\qquad
h_2(\alpha)=\frac{-47+45\sqrt{-7}}{8192},\qquad
A=\frac{1+64h_2}{1-64h_2}=\frac{9\sqrt{-7}}{35}.
\]
This together with Lemmas 4.1-4.2 yields
\[
\delta(\alpha)
=\frac{A}{\Omega^2}\Big(\theta\Omega-\frac{\Omega}{2\pi\,
\operatorname{Im}\alpha}\Big)
=\frac{9\sqrt{-7}}{35}\cdot\frac{8\sqrt{-7}}{63}
=\frac{72\cdot(-7)}{35\cdot63}
=-\frac{504}{2205}=-\frac{8}{35}.
\]
This proves the lemma.
\vskip0.2cm
\par{\bf Theorem 5.2} {\sl Let $p>7$ be a prime such that $p\e 1,2,4\mod 7$ and so $p=x^2+7y^2$. Then
$$\sum_{k=0}^{p-1}(35k+8)\f{\b{2k}k^2\b{4k}{2k}}{81^k}\e 8p\mod {p^3}\eqno{(5.3)}$$
and
$$\sum_{k=0}^{p-2}\f{\b{2k}k^2\b{4k}{2k}}{81^k(k+1)}\e -\f{100}{3}y^2+2p\pmod {p^3}.$$}
\par{\it Proof.} Putting $c=4y$, $d=x-y$, $e=1$, $\alpha=\f{1+\sqrt{-7}}4$, $N=2$ and $w=\chi=1$ in [2, Theorem 1.22] and then applying Lemma 4.3 and [20, Section 4] yields the first congruence. By [20, (1.7)],
$$\sum_{k=0}^{p-1}\f{\b{2k}k^2\b{4k}{2k}}{81^k}\e 4x^2-2p-\f{p^2}{4x^2}\mod {p^3}.\eqno{(5.4)}$$
This together with (5.3) gives
$$\sum_{k=0}^{p-1}k\f{\b{2k}k^2\b{4k}{2k}}{81^k}\e 
\f 8{35}\Big(3p-4x^2+\f{p^2}{4x^2}\Big)
\mod{p^3}.\eqno{(5.5)}$$
Now, putting $m=\f{81}{64}$ in Corollary 2.3 and then applying
(5.3), (5.5) and (3.3) yields
\begin{align*}
\sum_{k=0}^{p-2}\f{\b{2k}k^2\b{4k}{2k}}{81^k(k+1)}
&\e \sum_{k=0}^{p-1}\f{\b{2k}k^2\b{4k}{2k}}{81^k}
+\Big(\sum_{k=0}^{p-1}\f{\b{2k}k^2\b{4k}{2k}}{81^k}
\Big)^{-1}\Big(\f{175}{48}\big(\sum_{k=0}^{p-1}\f{k\b{2k}k^2
\b{4k}{2k}}{81^k}\big)^2+p^2\Big)
\\&\e 4x^2-2p-\f{p^2}{4x^2}+\Big(\f 1{4x^2}+\f p{8x^4}+\f{5p^2}{64x^6}\Big)\Big(\f 4{21}\big(3p-4x^2+\f{p^2}{4x^2}\big)^2+p^2\Big)
\\&\e -\f{100}{3}y^2+2p\pmod {p^3}.
\end{align*}
This proves the theorem.
\vskip0.2cm
\par In a similar way, one can prove the following results.
\skip0.2cm
\par{\bf Theorem 5.3} {\sl Let $p$ be a prime.
If $p=x^2+4y^2\e 1\pmod 4$, 
then$$\sum_{k=0}^{p-1}(7k+1)\f{\b{2k}k^2\b{4k}{2k}}{648^k}
\e p\mod{p^3}$$
and 
$$\sum_{k=0}^{p-2}\f{\b{2k}k^2\b{4k}{2k}}{648^k(k+1)}\e -\f {40}3y^2+2p\pmod {p^3}.$$}

\par{\bf Theorem 5.4} {\sl Let $p$ be a prime.
If $p=x^2+3y^2\e 1\pmod 3$, then
$$\sum_{k=0}^{p-1}(5k+1)\f{\b{2k}k^2\b{4k}{2k}}{(-144)^k}
\e \Ls p3p\mod {p^3}$$
and
$$\sum_{k=0}^{p-2}\f{\b{2k}k^2\b{4k}{2k}}{(-144)^k(k+1)}\e -16y^2+2p\pmod {p^3}.$$}

\par{\bf Theorem 5.5} {\sl Let $p$ be a prime.
If $p=x^2+2y^2\e 1,3\pmod 8$, then
$$\sum_{k=0}^{p-1}(40k+3)\f{\b{2k}k^2\b{4k}{2k}}{28^{4k}}\e 3\Ls p3p\mod{p^3}$$ and
$$\sum_{k=0}^{p-2}\f{\b{2k}k^2\b{4k}{2k}}{28^{4k}(k+1)}\e -284x^2+\big(142+144\sls p3\big)p\pmod {p^3}.$$}

\par{\bf Theorem 5.6} {\sl Let $p>7$ be a prime.
If $p=x^2+7y^2\e 1,2,4\pmod
7$, then
$$\sum_{k=0}^{p-1}(65k+8)\f{\b{2k}k^2\b{4k}{2k}} {(-3969)^k}
\e 8p\mod {p^3}$$
and
$$\sum_{k=0}^{p-2}\f{\b{2k}k^2\b{4k}{2k}} {(-3969)^k(k+1)}\e -\f{196}3y^2+2p\pmod {p^3}.$$}

\par{\bf Theorem 5.7} {\sl Let $p$ be a prime such that $p\e 1\mod 4$. Then
$$\sum_{k=0}^{p-1}(28k+3)\f{\b{2k}k^2\b{4k}{2k}}
{(-12288)^k}\e 3\Ls p3p\mod{p^3}$$ and
\begin{align*}&\sum_{k=0}^{p-2}\f{\b{2k}k^2\b{4k}{2k}}
{(-12288)^k(k+1)}
\\&\e\begin{cases} -144y^2+2p\pmod {p^3}
&\t{if $p\e 1\mod {12}$ and so
$p=x^2+9y^2$,}
\\72y^2-2p\pmod {p^3} &\t{if $p\e 5\mod{12}$ and so
$2p=x^2+9y^2$.}
\end{cases}\end{align*}}

\par{\bf Theorem 5.8} {\sl Let $p$ be a prime such that $\sls{-6}p=1$.
Then
$$\sum_{k=0}^{p-1}(8k+1)\f{\b{2k}k^2\b{4k}{2k}}
{48^{2k}}\e \Ls p3p\mod {p^3}$$ and
\begin{align*}\sum_{k=0}^{p-2}\f{\b{2k}k^2\b{4k}{2k}}
{48^{2k}(k+1)}\e\begin{cases}
-8y^2+2p\pmod {p^3}&\t{if $p=x^2+6y^2\e 1,7\pmod
{24}$,}\\4y^2-2p\pmod {p^3}&\t{if $p=2x^2+3y^2\e 5,11\pmod
{24}$}.\end{cases}
\end{align*}}

\par{\bf Theorem 5.9} {\sl Let $p$ be a prime such that $\sls{-10}p=1$. Then
$$\sum_{k=0}^{p-1}(10k+1)\f{\b{2k}k^2\b{4k}{2k}}{12^{4k}}
\e \Ls{-2}pp\mod {p^3}$$ and
\begin{align*}&\sum_{k=0}^{p-2}\f{\b{2k}k^2\b{4k}{2k}}{12^{4k}(k+1)}
\\&\e\begin{cases}\f{392}{3}y^2+2p-p^2\pmod {p^3}&\t{if $p=x^2+10y^2\e
1,9,11,19\pmod {40}$,}
\\-\f{196}{3}y^2-2p-p^2\pmod {p^3}&\t{if $p=2x^2+5y^2\e 7,13,23,37\pmod {40}$.}
\end{cases}\end{align*}}

\par{\bf Theorem 5.10} {\sl Let $p>3$ be a prime and $\sls{-13}p=1$.
Then
$$\sum_{k=0}^{p-1}(260k+23)\f{\b{2k}k^2\b{4k}{2k}}{(-82944)^k}
\e 23\Ls{-1}pp\mod {p^3}$$
and
\begin{align*}&\sum_{k=0}^{p-2}\f{\b{2k}k^2\b{4k}{2k}}{(-82944)^k(k+1)}
\\&\e \begin{cases}
-\f{2272}3y^2+2p\pmod {p^3}&\t{if $\sls{-1}p=\sls{13}p=1$ and so $p=x^2+13y^2$,}
\\\f{1136}3y^2-2p\pmod {p^3}&\t{if $\sls{-1}p=\sls{13}p=-1$ and so
$2p=x^2+13y^2$.}\end{cases}\end{align*}}

\par{\bf Theorem 5.11} {\sl Let $p$ be a prime with $p\not=2,3,7$ and $\sls{-37}p=1$. Then
$$\sum_{k=0}^{p-1}(21460k+1123)
\f{\b{2k}k^2\b{4k}{2k}}{(-199148544)^k}\e 1123p\mod {p^3}$$
and
\begin{align*}&\sum_{k=0}^{p-2}\f{\b{2k}k^2\b{4k}{2k}}{(-199148544)^k(k+1)}
\\&\e \begin{cases}
-\f{5044960}3y^2+2p\pmod {p^3}&\t{if $\sls{-1}p=\sls{37}p=1$ and so $p=x^2+37y^2$,}
\\ \f{2522480}3y^2-2p\pmod {p^3}&\t{if $\sls{-1}p=\sls{37}p=-1$ and so
$2p=x^2+37y^2$.}
\end{cases}\end{align*}}

\par{\bf Theorem 5.12} {\sl Let $p$ be a prime with $p\not=2,3,11$ and $\sls{-22}p=1$. Then
$$\sum_{k=0}^{p-1}(280k+19)\f{\b{2k}k^2\b{4k}{2k}}{1584^{2k}}
\e 19\Ls p{11}p\mod {p^3}$$
and
\begin{align*}&\sum_{k=0}^{p-2}\f{\b{2k}k^2\b{4k}{2k}}{1584^{2k}(k+1)}
\\&\e \begin{cases}\f{63272}3y^2+2p\pmod {p^3}&\t{if $\ls{2}p=\sls{-11}p=1$ and so $p=x^2+22y^2$,}
\\ -\f{31636}3y^2-2p\pmod {p^3}&\t{if $\ls{2}p=\sls{-11}p=-1$ and so $p=2x^2+11y^2$.}
\end{cases}\end{align*}}

\par{\bf Theorem 5.13} {\sl
Let $p$ be a prime, $p\not=2,3,11$ and $\ls{-58}p=1$.
Then
$$\sum_{k=0}^{p-1}(26390k+1103)\f{\b{2k}k^2\b{4k}{2k}}
{396^{4k}}\e 1103\Ls{-2}pp\mod {p^3}$$
and
\begin{align*}
&\sum_{k=0}^{p-2}\f{\b{2k}k^2\b{4k}{2k}}{396^{4k}(k+1)}
\\&\e \begin{cases}
\f{622903112}3y^2+2p\pmod {p^3}&\t{if $\ls{-2}p=\sls{29}p=1$ and so $p=x^2+58y^2$,}
\\ -\f{311451556}3y^2-2p\pmod {p^3}&\t{if $\sls{-2}p=\sls{29}p=-1$ and so
$p=2x^2+29y^2$.}\end{cases}\end{align*}
}

\par{\bf Theorem 5.14} {\sl
Let $p>5$ be a prime of the form $4k+1$. Then
$$\sum_{k=0}^{p-1}(644k+41)
\f{\b{2k}k^2\b{4k}{2k}}{(-6635520)^k}\e 41\Ls{-5}pp\mod{p^3}$$
and
\begin{align*}
&\sum_{k=0}^{p-2}\f{\b{2k}k^2\b{4k}{2k}}{(-6635520)^k(k+1)}
\\&\e \begin{cases}
-\f{168400}3y^2+2p\pmod {p^3}&\t{if 
 $p\e 1,9\pmod {20}$
and so $p=x^2+25y^2$}, 
\\ \f{84200}3y^2-2p\pmod {p^3}&\t{if $p\e 13,17\pmod {20}$ and so $2p=x^2+25y^2$.}
\end{cases}\end{align*}}
\par{\bf Remark 5.1} Most of the Van Hamme-type congruences in Theorems 5.1-5.14 were conjectured by Z.W. Sun in [21] and [23]. Combining Theorem 2.2 with Theorems 5.1-5.14 one may deduce the congruences on $$\sum_{k=0}^{p-1}w(k)\f{\b{2k}k^2\b{4k}{2k}}{m^k}\q
(w(k)\in\{k^2,k^3,2k-1\})$$ modulo $p^3$  and  $\sum_{k=0}^{p-2}\b{2k}k^2\b{4k}{2k}/(k+1)^r$ modulo $p^3$ for $r=2,3$, where $p>3$ is a prime and $m\in\Bbb Z$ with $p\nmid m$.

%\par{\bf Declaration of competing interest}
%\par The author declares no conflicts of interest regarding the publication of this paper.
\vskip0.2cm
\par{\bf Funding} The author was supported by the
National Natural Science Foundation of China (Grant No. 12271200).
\vskip0.2cm
%\par{\bf Data Availability}
%\par No data was used for the research described in the article.


\begin{thebibliography}{99}

\bibitem [1]{} B.C. Berndt, R.J. Evans and K.S. Williams, {\it
Gauss and Jacobi Sums}, Wiley, New York, 1998.

\bibitem [2] {} F. Beukers, {\it Supercongruences using modular forms}, Algebra \& Number Theory {\bf 20}\,(2026), 
    1727-1768.
    
    \bibitem [3] {} D.A. Cox, {\it Primes of the form $x^2+ny^2$: Fermat, Class Field Theory, and Complex Multiplication, 2nd edition}, in: Pure and Applied Mathematics, Wiley, New York, 2013.

 \bibitem [4] {} J. Guillera and W. Zudilin, {\it ``Divergent" Ramanujan-type
supercongruences}, Proc. Amer. Math. Soc. {\bf 140}(2012), 765-777.

\bibitem [5] {}  L. Long, {\it Hypergeometric evaluation identities and
supercongruences}, Pacific J. Math. {\bf 249}(2011), 405-418.

\bibitem [6] {} G.-S. Mao, {\it On some congruence conjectures involving binary quadratic forms}, Publ. RIMS Kyoto Univ. {\bf 61}(2025), 611-625.

\bibitem [7] {} E. Mortenson, {\it A p-adic supercongruence conjecture of
 van Hamme}, Proc. Amer. Math. Soc. {\bf 136}(2008), 4321-4328.

\bibitem [8] {} R. Osburn and W. Zudilin,
{\it On the (K.2) supercongruence of Van Hamme}, J. Math. Anal.
Appl. {\bf 433}(2016), 706-711.

\bibitem [9]{}
L.\,J. Slater, \emph{Generalized hypergeometric functions},
Cambridge Univ. Press, Cambridge, 1966.

\bibitem [10] {}
Z.H. Sun, {\it Congruences concerning Legendre polynomials II}, J.
Number Theory {\bf 133}(2013), 1950-1976.

\bibitem [11] {}
Z.H. Sun, {\it Congruences involving $\b{2k}k^2\b{3k}k$}, J. Number
Theory {\bf 133}(2013), 1572-1595.

\bibitem [12] {}
Z.H. Sun, {\it Legendre polynomials and supercongruences}, Acta
Arith. {\bf 159}(2013), 169-200.

\bibitem [13] {}
Z.H. Sun, {\it Generalized Legendre polynomials and related
supercongruences}, J. Number Theory {\bf 143}(2014), 293-319.

\bibitem [14] {} Z.H. Sun, {\it Congruences involving binomial
coefficients and Ap\'ery-like numbers}, Publ. Math. Debrecen {\bf
96}(2020), 315-346.



\bibitem [15] {} Z.H. Sun, {\it
New conjectures involving binomial coefficients and Ap\'ery-like numbers}, arXiv:2111.04538.

\bibitem [16] {} Z.H. Sun, {\it Supercongruences involving Ap\'ery-like numbers and binomial
coefficients}, AIMS Math. {\bf 7}(2022), 2729-2781.

\bibitem [17] {} Z.H. Sun, {\it Supercongruences involving products of three binomial coefficients}, Rev. Real Acad. Cienc. Exactas Fis. Nat. Ser. A-Mat. {\bf 117}(2023), Art.131, 26pp.

\bibitem [18] {} Z.H. Sun, {\it Binomial Coefficients, Recurrence Sequences and Congruences} (Chinese), Science Press, Beijing, 2025.

    \bibitem [19] {} Z.H. Sun, {\it Congruences involving the products of three binomial coefficients,} Res. Number Theory {\bf 12}\,(2026), Art. 68, 37pp.


    \bibitem [20] {}
    Z.H. Sun and D. Ye, {\it Supercongruences involving products of three binomial coefficients via Beukers' method}, Results in Math. {\bf 80}(2025), Art.143, 26pp.

\bibitem [21] {} Z.W. Sun, {\it Super congruences and Euler
numbers}, Sci. China Math. {\bf 54}(2011), 2509-2535.

\bibitem [22] {} Z.W. Sun, {\it On sums involving products of three
binomial coefficients}, Acta Arith. {\bf 156}(2012), 123-141.

\bibitem  [23] {} Z.W. Sun, {\it Conjectures and results on $x^2$
mod $p^2$ with $4p=x^2+dy^2$}, in: Number Theory and Related Areas
(eds., Y. Ouyang, C. Xing, F. Xu and P. Zhang), Higher Education
Press $\&$ International Press, Beijing and Boston, 2013,
pp.149-197.


\end{thebibliography}
\end{document}